\documentclass[preprint,12pt,authoryear]{elsarticle}

\usepackage{amsmath}
\usepackage{amssymb}
\usepackage{amsthm}
\usepackage{mathtools}
\usepackage{xcolor}
\usepackage{graphicx}
\usepackage{float}
\journal{}

\newtheorem{assumption}{Assumption}
\newtheorem{lemma}{Lemma}
\newtheorem{proposition}{Proposition}
\newtheorem{theorem}{Theorem}

\newcounter{modelalgorithm}
\newenvironment{algorithmblock}[1]
{%
\par\medskip
\noindent\begin{minipage}{\linewidth}
\refstepcounter{modelalgorithm}%
\noindent\hrule
\medskip
\noindent\textbf{Algorithm \themodelalgorithm. #1}
\par\medskip
}
{%
\medskip
\noindent\hrule
\end{minipage}
\par\medskip
}

\begin{document}

\begin{frontmatter}

\title{{Reinforcement learning for constrained stochastic linear--quadratic control in continuous time}}

\author[1]{Weiping Wu}
\author[1]{Qitong Pan}
\author[2]{Xun Li$^*$}
\author[1]{Peng Wu$^*$}

\affiliation[1]{organization={School of Economics and Management, Fuzhou University}, city={Fuzhou}, postcode={350108}, state={Fujian}, country={China}}

\affiliation[2]{organization={Department of Applied Mathematics, The Hong Kong Polytechnic University}, city={Hong Kong},  country={China}}

\begin{abstract}
{
We study continuous-time stochastic linear--quadratic control with a scalar state and nonnegative
control constraints. With quadratic running costs, randomized exploration can change the optimal
mean feedback. We introduce an auxiliary problem without entropy in which the conditional policy
mean is nonnegative and the covariance is prescribed in proportion to the squared state. Moment
reduction and state separation yield two regional Riccati equations and characterize the
finite-horizon value function and optimal conditional mean. We give a sufficient control--cost
compatibility condition under which this mean coincides with the optimal feedback of the original
constrained problem. For time-invariant coefficients, we characterize the undiscounted
infinite-horizon optimum and establish mean-square asymptotic stability under the stated
assumptions. For unknown dynamics, we develop model-free actor--critic schemes for finite and
infinite horizons. Numerical experiments compare the learned policies with their analytical targets
and examine the effect of exploration on the optimal mean.
}
\end{abstract}


\end{frontmatter}

\section{Introduction}
\label{sec:introduction}

{For classical unconstrained stochastic linear--quadratic (LQ) problems, Riccati equations provide a means of
constructing optimal controls in linear state-feedback form; see \citet[Chapter~6]{YongZhou1999}.} Its connection
with mean--variance portfolio selection was established by \citet{ZhouLi2000}, using the indefinite
stochastic LQ theory of \citet{ChenLiZhou1998}. Control constraints introduce difficulties even when
the model coefficients are known. In the no-shorting problem, \citet{LiZhouLim2002} showed that the
usual completion-of-squares argument does not apply directly and that the associated Hamilton--Jacobi--Bellman
(HJB) equation need not have a smooth solution. They characterized the solution through two Riccati
equations and viscosity analysis. \citet{HuZhou2005} developed cone-constrained stochastic LQ control
with random coefficients through two extended stochastic Riccati equations, while \citet{LiXu2016}
studied wealth and portfolio constraints by a martingale approach. For scalar-state systems with
linear control constraints, \citet{WuEtAl2020} established state separation and derived regional
Riccati equations, including an infinite-horizon characterization under additional conditions.

{
These model-based approaches rely on specified system dynamics \citep{WuEtAl2020,WangZhou2020}.
Their use in mean--variance portfolio selection faces two distinct difficulties. First, financial
model specification must balance applicability and tractability, and models that retain more
realistic market features can become difficult to solve with classical stochastic control tools
\citep{HamblyXuYang2023}. 
Second, implementing model-based portfolio policies requires coefficient estimates. Accurate estimation is difficult,
and portfolio allocations are sensitive to estimation errors \citep[Section~2.2]{WangZhou2020}.
Model-free reinforcement learning (RL) can bypass explicit coefficient estimation and learn controls
from observed interactions \citep{HamblyXuYang2023,WangZhou2020}. However, the squared mean in a
variance objective makes it nonlinear in expectation. The usual learning rules for expected
cumulative rewards therefore cannot be applied directly to mean--variance criteria
\citep[Introduction]{WangZhou2020}.}

{
Exploratory stochastic control treats exploration explicitly in continuous time by replacing
actions with probability distributions \citep{WangZariphopoulouZhou2020}. The LQ structure provides
analytical optimal exploratory policies in the settings of \citet{WangZariphopoulouZhou2020,WangZhou2020}.
Auxiliary formulations without entropy have also recovered original optimal controls through
policy means for specific investment criteria \citep{DaiEtAl2026,LiWangZhou2026}. We study this
auxiliary approach for scalar-state stochastic LQ control with nonnegative control constraints and
quadratic running costs. We characterize its finite-horizon value and optimal conditional mean.
Quadratic running costs can make exploration change the optimal mean, so its relationship with the
original constrained feedback requires a separate analysis (Theorems~\ref{thm:main}
and~\ref{thm:alignment_consistency}). We also study the corresponding undiscounted infinite-horizon
problem and construct learning methods for unknown dynamics
(Sections~\ref{subsec:infinite_horizon} and~\ref{sec:unknown_dynamics}).}

{
The main contributions of this paper are as follows.
\begin{enumerate}
\item \textbf{Auxiliary formulation and finite-horizon solution.}
We formulate an auxiliary problem without entropy for scalar-state stochastic LQ control with
{nonnegative controls, using general randomized feedback policies} with nonnegative means and prescribed
state-dependent covariance. Moment reduction and state separation characterize its finite-horizon
value function and optimal conditional mean
(Section~\ref{subsec:randomized_policies} and Theorem~\ref{thm:main}).
\item \textbf{Consistency with the original constrained problem.}
We give a sufficient control--cost compatibility condition under which the optimal conditional
mean of the auxiliary problem coincides with the original optimal constrained feedback
(Theorem~\ref{thm:alignment_consistency}).
\item \textbf{Undiscounted infinite-horizon control.}
We characterize the infinite-horizon optimal value and conditional mean, and establish mean-square
asymptotic stability under the stated assumptions (Theorem~\ref{thm:infinite_horizon}).
\item \textbf{Learning with unknown dynamics.}
Using the martingale methods of \citet{JiaZhou2022PE,JiaZhou2022PG}, we develop finite- and
infinite-horizon model-free actor--critic schemes. Numerical experiments compare the learned
policies with their analytical targets
(Sections~\ref{sec:unknown_dynamics} and~\ref{sec:numerical_experiments}).
\end{enumerate}
}

{
The paper is organized as follows. Section~\ref{subsec:related_literature} reviews the related
literature. Section~\ref{sec:problem_formulation} defines the original and auxiliary problems.
Section~\ref{sec:main_results} establishes the exploratory-control results.
Section~\ref{sec:unknown_dynamics} presents the learning methods, and
Section~\ref{sec:numerical_experiments} reports the numerical studies.
Section~\ref{sec:conclusion} concludes the paper. The appendices provide the proofs and learning
derivations.}

{
\subsection{Related literature}
\label{subsec:related_literature}}

{
Mean--variance (MV) portfolio selection minimizes the variance of terminal wealth for a
prescribed expected terminal wealth \citep{WangZhou2020}. \citet{ZhouLi2000} embedded the
continuous-time MV portfolio problem into a class of auxiliary stochastic LQ problems. In this
formulation, wealth is the scalar state and portfolio positions are the controls
\citep[Section~2.1]{WangZhou2020}.
\citet{HamblyXuYang2023} reviewed applications of RL to portfolio optimization and hedging.
Continuous-time portfolio learning without explicit coefficient estimation has been developed by
\citet{WangZhou2020,HuangJiaZhou2024}. Variance-related
objectives and constraints have also been studied for cumulative rewards
\citep[Section~3]{PrashanthGhavamzadeh2016}.
\citet{TamarDiCastroMannor2016} developed methods for evaluating variance under a fixed policy.
\citet{TamarDiCastroMannor2012,PrashanthGhavamzadeh2016} studied variance-related objectives and
constraints using local policy optimization. These evaluation and local optimization results do
not provide a general guarantee of a globally optimal mean--variance policy.}

{
To model the exploration--exploitation trade-off explicitly in continuous time,
\citet{WangZariphopoulouZhou2020} introduced an exploratory stochastic control formulation.
In this formulation, the drift, diffusion variance, and instantaneous reward are averaged with
respect to action distributions.
Entropy regularization balances exploration and exploitation. In the unconstrained LQ setting of
\citet{WangZariphopoulouZhou2020}, the optimal policy is Gaussian. Its mean represents exploitation,
and its variance represents exploration. The policy and value function converge to their classical
counterparts as the exploration weight vanishes. For finite-horizon mean--variance portfolio
selection, \citet{WangZhou2020} used the LQ structure to derive a globally optimal exploratory policy.
They also proved a policy improvement theorem. Under exact policy evaluation and their stated
Gaussian initialization, policy iteration converges to that optimum. These results give a global
characterization for this specific LQ problem beyond the local optimization results discussed above.
Martingale methods subsequently provided continuous-time policy evaluation, policy-gradient
representations, and $q$-learning \citep{JiaZhou2022PE,JiaZhou2022PG,JiaZhou2023Q}.
Related analyses address exploratory HJB equations \citep{TangZhangZhou2022}, entropy annealing
\citep{SethiSiskaZhang2025}, LQ policy convergence \citep{GiegrichEtAl2024,GuoLiXu2026}, and the
accuracy of discretely sampled stochastic policies \citep{JiaOuyangZhang2026}.}

{
Control constraints also complicate exploration. Constraining the support of randomized policies
can exclude ordinary Gaussian policies and make theoretical analysis more difficult
\citep{LiWangZhou2026}. The interval constraints studied by \citet{ChauNguyenNguyen2026} lead to
truncated Gaussian policies with additional normalization and boundary terms. An alternative is
to formulate an auxiliary problem without entropy. \citet{DaiEtAl2026} proved that the mean of an
optimal randomized policy solves their original unconstrained Merton problem. For no-shorting
mean--variance portfolio selection, \citet{LiWangZhou2026} imposed the constraint on the policy
mean and proved that, with a suitable prescribed covariance, the optimal mean recovers the original
optimal portfolio. These investment results do not establish the mean-policy relationship with
quadratic running costs or characterize the undiscounted infinite-horizon auxiliary problem with
prescribed state-dependent covariance for constrained stochastic LQ control.}

\section{Problem Formulation}
\label{sec:problem_formulation}

{We first formulate the {original stochastic LQ problem with nonnegative controls} and then define its
randomized auxiliary problem. In the auxiliary problem, the conditional action mean satisfies the
constraint and a prescribed covariance governs exploration. }

\subsection{Notation}
\label{subsec:notation}

Frequently used notation in this paper includes:
\begin{itemize}
\item $M^\top$: the transpose of a matrix or vector $M$;
\item {$\|M\|:=\left(\sum_{i,j}m_{ij}^2\right)^{1/2}$ for a matrix or vector
$M=(m_{ij})$;}
\item $\mathbb R$, $\mathbb R^n$, and $\mathbb R^{n\times m}$: the sets of real numbers, real $n$-vectors, and real $n\times m$ matrices, respectively;
\item $\mathbb R_+:=[0,\infty)$ and $\mathbb R_+^n:=[0,\infty)^n$;
\item $\mathbb S^n:=\{M\in\mathbb R^{n\times n}:M^\top=M\}$, with $\mathbb S_+^n$ and $\mathbb S_{++}^n$ denoting its positive semidefinite and positive definite subsets;
\item $M\succeq N$ and $M\succ N$: $M-N$ is positive semidefinite and positive definite, respectively;
\item $a\ge0$: a componentwise inequality when $a$ is a vector;
\item $\mathcal N(\mu,\Sigma)$: the Gaussian probability distribution with mean $\mu$ and covariance matrix $\Sigma$;
\item {$\mathbf 1_E$: the indicator of an event $E$;}
\item $\operatorname{tr}(M)$: the trace of a square matrix $M$;
\item $I_n$: the $n\times n$ identity matrix;
\end{itemize}

Let $T>0$ be a finite horizon. Consider a complete filtered probability space $(\Omega,\mathcal F,\mathbb
F,\mathbb P)$ satisfying the usual conditions and carrying an $m$-dimensional standard Brownian motion
$W=(W_1,\ldots,W_m)^\top$, where $\mathbb F=\{\mathcal F_t\}_{0\le t\le T}$.  Here $\mathcal F_t$ is the
augmented $\sigma$-algebra generated by the Brownian history up to time $t$ and represents the
information available at time $t$. {The symbol $\mathbb E$ denotes expectation.}

Given a Hilbert space $\mathcal H$ with norm $\|\cdot\|_{\mathcal H}$, define the Banach space
\begin{equation}
L_{\mathbb F}^2(0,T;\mathcal H):=
\left\{\varphi:\ \begin{array}{l}
\varphi\text{ is an }\mathbb F\text{-adapted, }\mathcal H\text{-valued measurable process},\\
\displaystyle \mathbb E\int_0^T\|\varphi(t)\|_{\mathcal H}^2dt<\infty
\end{array}\right\},
\label{eq:L2_adapted_space}
\end{equation}
equipped with the norm
\[
\|\varphi\|_{L_{\mathbb F}^2(0,T;\mathcal H)}
:=\left(\mathbb E\int_0^T\|\varphi(t)\|_{\mathcal H}^2dt\right)^{1/2}.
\]
For $t\in[0,T]$, the space $L_{\mathbb F}^2(t,T;\mathcal H)$ is defined analogously. {The nonnegative control
constraint is incorporated} through the subset
\begin{equation}
L_{\mathbb F}^2(t,T;\mathbb R_+^n)
:=\left\{u\in L_{\mathbb F}^2(t,T;\mathbb R^n):
{u(s)\in\mathbb R_+^n,\ s\in[t,T]}\right\}.
\label{eq:L2_noshorting_space}
\end{equation}

\subsection{{Original stochastic LQ problem with nonnegative controls}}
\label{subsec:original_problem}

We first state the finite-horizon problem from the initial time $0$. For a
given initial state $x_0\in\mathbb R$, consider
\begin{equation}
\left\{
\begin{aligned}
dX(t)={}&\left[A(t)X(t)+B(t)^\top u(t)\right]dt\\
&+\left[X(t)C(t)^\top+u(t)^\top D(t)\right]dW(t), \qquad t\in[0,T],\\
X(0)={}&x_0,
\end{aligned}
\right.
\label{eq:classical_state}
\end{equation}
where $A(t)\in\mathbb R$, $B(t)\in\mathbb R^n$, $C(t)\in\mathbb R^m$, $D(t)\in\mathbb R^{n\times m}$, and
$u\in L_{\mathbb F}^2(0,T;\mathbb R_+^n)$ is the control.

Let $R(t)\in\mathbb S^n$, $S(t)\in\mathbb R^n$, $q(t)\in\mathbb R$,
and $q_T\in\mathbb R$. The running cost is
\begin{equation}
\ell(t,x,u):=u^\top R(t)u+2xS(t)^\top u+q(t)x^2.
\label{eq:classical_running_cost}
\end{equation}
The superscript $0$ denotes the original problem without randomized exploration.
The performance criterion from time $0$ is
\begin{equation}
J^0(x_0;u):=\mathbb E\left[\int_0^T\ell\bigl(t,X(t),u(t)\bigr)dt+q_TX(T)^2\right].
\label{eq:classical_performance}
\end{equation}
The original finite-horizon constrained stochastic LQ problem is therefore
\begin{equation}
V^0(0,x_0):=\inf_{u\in L_{\mathbb F}^2(0,T;\mathbb R_+^n)}J^0(x_0;u).
\label{eq:classical_value}
\end{equation}

\begin{assumption}
{The coefficient functions $A(\cdot)$, $B(\cdot)$, $C(\cdot)$, $D(\cdot)$, $R(\cdot)$, $S(\cdot)$, and
$q(\cdot)$ are deterministic, measurable, and bounded on $[0,T]$.}
\end{assumption}

\begin{assumption}
For every $t\in[0,T]$,
\begin{equation}
Q(t):=\begin{pmatrix}R(t)&S(t)\\S(t)^\top&q(t)\end{pmatrix}\succeq0, \qquad q_T\ge0.
\label{eq:cost_convexity}
\end{equation}
\end{assumption}

{
\begin{assumption}
\label{ass:control_weight_positive_definite}
{The control-weight matrix is pointwise positive definite on $[0,T]$, meaning that every nonzero control
direction has a strictly positive quadratic weight:}
\begin{equation}
R(t)\succ0, \qquad t\in[0,T].
\label{eq:control_weight_positive_definite}
\end{equation}
\end{assumption}
}

\begin{assumption}
There exists $\delta>0$ such that
\begin{equation}
D(t)D(t)^\top\succeq\delta I_n, \qquad t\in[0,T].
\label{eq:diffusion_nondegeneracy}
\end{equation}
\end{assumption}

The boundedness in Assumption~1 provides the uniform coefficient control needed for the usual
square-integrability estimates of the state process.

For dynamic programming, we next introduce the continuation problem at an arbitrary state--time pair
$(t,x)\in[0,T]\times\mathbb R$. {For $u\in L_{\mathbb F}^2(t,T;\mathbb R_+^n)$, let $X$ solve}
\[
\left\{
\begin{aligned}
dX(s)={}&[A(s)X(s)+B(s)^\top u(s)]ds\\
&+[X(s)C(s)^\top+u(s)^\top D(s)]dW(s),\qquad s\in[t,T],\\
X(t)={}&x.
\end{aligned}
\right.
\]
For a controlled state initialized by $X(t)=x$, $\mathbb E_{t,x}$ denotes expectation conditional on this
initialization. The corresponding cost-to-go and value function are
\[
\begin{aligned}
J(t,x;u)&:=\mathbb E_{t,x}\left[\int_t^T\ell\bigl(s,X(s),u(s)\bigr)ds+q_TX(T)^2\right],\\
V^0(t,x)&:=\inf_{u\in L_{\mathbb F}^2(t,T;\mathbb R_+^n)}J(t,x;u).
\end{aligned}
\]
{Thus $V^0(0,x_0)$ is the value of the original problem. {This model with nonnegative controls belongs} to the
scalar-state constrained stochastic LQ class of \citet{WuEtAl2020}. Under Assumptions~1--4, its optimal
feedback is characterized by two regional quadratic programs and Riccati equations. We use this
feedback below to assess the effect of exploration.}

{We next introduce a randomized auxiliary problem whose conditional action mean satisfies the
{same nonnegative control constraint.}}

\subsection{General randomized feedback policies}
\label{subsec:randomized_policies}

{The conditional mean of a randomized policy carries the control constraint, while its
covariance specifies the exploration level.}

{Following \citet[Definition~6.4]{Villani2009}, let $\mathcal B(\mathbb R^n)$ denote the Borel
$\sigma$-algebra of $\mathbb R^n$, and define}
{
\[
\mathcal P_2(\mathbb R^n):=
\left\{\nu:\ \begin{array}{l}
\nu\text{ is a probability measure on }(\mathbb R^n,\mathcal B(\mathbb R^n)),\\[1mm]
\displaystyle\int_{\mathbb R^n}\|u\|^2\nu(du)<\infty
\end{array}\right\}.
\]
}
{Thus $\mathcal P_2(\mathbb R^n)$ is the set of Borel probability measures on $\mathbb R^n$ with finite
second moment. The randomized action law need not admit a density, while its finite second moment makes
the mean, covariance matrix, and quadratic terms well defined. {A randomized feedback policy is a
measurable map}}
\begin{equation}
\pi:[0,T]\times\mathbb R\longrightarrow\mathcal P_2(\mathbb R^n), \qquad (s,x)\longmapsto\pi_s(\cdot\mid x).
\label{eq:randomized_feedback_policy}
\end{equation}
For a policy $\pi$, define its conditional mean and {covariance matrix} by
\begin{equation}
a^\pi(s,x):=\int_{\mathbb R^n}u\,\pi_s(du\mid x),
\label{eq:policy_mean}
\end{equation}
\begin{equation}
\Sigma^\pi(s,x):=\int_{\mathbb R^n}\bigl(u-a^\pi(s,x)\bigr)\bigl(u-a^\pi(s,x)\bigr)^\top\pi_s(du\mid x).
\label{eq:policy_covariance}
\end{equation}
Hence,
\begin{equation}
\int_{\mathbb R^n}uu^\top\pi_s(du\mid x)=a^\pi(s,x)a^\pi(s,x)^\top+\Sigma^\pi(s,x).
\label{eq:policy_second_moment}
\end{equation}

{The nonnegative control constraint is imposed} on the conditional mean,\footnote{{The relationship between this
auxiliary mean constraint and the {original nonnegative control constraint} will be established below. In
particular, under the control--cost compatibility condition, Theorem~\ref{thm:alignment_consistency}
shows that the optimal conditional mean of the auxiliary problem coincides with the optimal
feedback control of the original problem.}}
\begin{equation}
a^\pi(s,x)\in\mathbb R_+^n,
\label{eq:mean_noshorting}
\end{equation}
{The original control must be componentwise nonnegative. In the auxiliary problem, this
constraint applies to $a^\pi(s,x)$ rather than to each sampled action. Samples may therefore have
negative components during exploration, whereas the extracted mean control belongs to
$\mathbb R_+^n$.}

We prescribe the exploration covariance as
\begin{equation}
\Sigma^\pi(s,x)=x^2\Psi(s),
\label{eq:prescribed_covariance}
\end{equation}
where $\Psi:[0,T]\to\mathbb S_+^n$ is deterministic, measurable, and bounded on $[0,T]$, i.e.,
$\sup_{t\in[0,T]}\|\Psi(t)\|<\infty$. This boundedness keeps the prescribed exploratory second moments
finite and preserves the quadratic-growth structure of the exploratory diffusion and running cost.

{To construct the state dynamics under a randomized policy, we average the drift and the quadratic
variation generated by the coefficients below. For $u\in\mathbb R^n$, let $C_j(s)$ be the $j$th component
of $C(s)$ and $D_{\cdot j}(s)$ the $j$th column of $D(s)$. Define}
\begin{equation}
b(s,x,u):=A(s)x+B(s)^\top u,
\label{eq:drift_coefficient}
\end{equation}
\begin{equation}
\sigma_j(s,x,u):=xC_j(s)+D_{\cdot j}(s)^\top u, \qquad j=1,\ldots,m.
\label{eq:diffusion_coefficient}
\end{equation}
{Under a randomized policy $\pi$, define the coefficients of the exploratory dynamics by}
\begin{equation}
\widetilde b^\pi(s,x):=\int_{\mathbb R^n}b(s,x,u)\,\pi_s(du\mid x),
\label{eq:exploratory_drift}
\end{equation}
\begin{equation}
\bigl(\widetilde\sigma_j^\pi(s,x)\bigr)^2:=\int_{\mathbb R^n}\sigma_j^2(s,x,u)\,\pi_s(du\mid x), \qquad j=1,\ldots,m.
\label{eq:exploratory_diffusion}
\end{equation}
{
\begin{proposition}
\label{prop:exploratory_characteristics}
Let $\pi$ be a randomized policy satisfying
\eqref{eq:prescribed_covariance}, and write $a=a^\pi(s,x)$. Then
\begin{equation}
\widetilde b^\pi(s,x)=A(s)x+B(s)^\top a,
\label{eq:drift_reduction}
\end{equation}
and
\begin{equation}
\sum_{j=1}^m\bigl(\widetilde\sigma_j^\pi(s,x)\bigr)^2
=\|xC(s)+D(s)^\top a\|^2
+x^2\operatorname{tr}\!\left[D(s)D(s)^\top\Psi(s)\right].
\label{eq:variance_reduction}
\end{equation}
\end{proposition}

\begin{proof}
See \ref{app:exploratory_characteristics}.
\end{proof}
}
{The exploratory dynamics under $\pi$ are}
\begin{equation}
\begin{cases}
dX^\pi(s)=\widetilde b^\pi\bigl(s,X^\pi(s)\bigr)ds+\displaystyle\sum_{j=1}^m\widetilde\sigma_j^\pi\bigl(s,X^\pi(s)\bigr)dW_j(s),\\
X^\pi(t)=x,
\end{cases}
\qquad s\in[t,T].
\label{eq:exploratory_state}
\end{equation}
{Here $\widetilde\sigma_j^\pi$ is the nonnegative square root specified in
\eqref{eq:exploratory_diffusion}.}

{The running cost under $\pi$ is}
\begin{equation}
\overline\ell^\pi(s,x):=\int_{\mathbb R^n}\ell(s,x,u)\,\pi_s(du\mid x).
\label{eq:running_cost_under_policy}
\end{equation}
The auxiliary exploratory problem is
\begin{equation}
J^\Psi(t,x;\pi):=\mathbb E_{t,x}\left[\int_t^T\overline\ell^\pi\bigl(s,X^\pi(s)\bigr)ds+q_TX^\pi(T)^2\right],
\label{eq:exploratory_performance}
\end{equation}
\begin{equation}
V^\Psi(t,x):=\inf_{\pi\in\Pi_{\Psi,+}[t,T]}J^\Psi(t,x;\pi).
\label{eq:exploratory_problem}
\end{equation}
{Here $\Pi_{\Psi,+}[t,T]$ denotes the set of {measurable randomized feedback policies} defined above that satisfy
\eqref{eq:mean_noshorting}--\eqref{eq:prescribed_covariance}, generate a mean process
$a^\pi(\cdot,X^\pi(\cdot))\in L_{\mathbb F}^2(t,T;\mathbb R_+^n)$, and make \eqref{eq:exploratory_state}
well posed and the criterion \eqref{eq:exploratory_performance} finite.}

{The admissible policies are not restricted to a parametric family. The next section shows
that the LQ structure and prescribed covariance reduce their optimization to an optimization over the
conditional mean.}

\section{Main Results}
\label{sec:main_results}

{We first reduce the auxiliary problem to conditional-mean optimization and characterize
its finite-horizon value by state separation. We then identify a condition for equality between the
exploratory optimal mean and the original optimal feedback. Finally, for constant coefficients, we
derive the {infinite-horizon optimal feedback with constant regional gains}. These results provide the targets for the learning
methods in Section~\ref{sec:unknown_dynamics}.}

\subsection{From distributional optimization to mean optimization}
\label{subsec:distribution_to_mean}

{For $\nu\in\mathcal P_2(\mathbb R^n)$, define its mean vector and covariance matrix by}
{
\begin{equation}
m_\nu:=\int_{\mathbb R^n}u\,\nu(du),
\label{eq:law_mean}
\end{equation}
\begin{equation}
\operatorname{Cov}_\nu(u)
:=\int_{\mathbb R^n}
\left(u-m_\nu\right)\left(u-m_\nu\right)^\top\nu(du).
\label{eq:law_covariance}
\end{equation}
}
{At a fixed state--time pair $(t,x)$, the measure $\nu$ represents
a candidate conditional action law. Define}
\begin{equation}
\mathcal P_{\Psi,+}(t,x):=\left\{\nu\in\mathcal P_2(\mathbb R^n):\ {m_\nu\in\mathbb R_+^n},\quad \operatorname{Cov}_\nu(u)=x^2\Psi(t)\right\}.
\label{eq:local_policy_set}
\end{equation}

{{For a candidate value function $v$}\footnote{The superscript $0$ denotes the original constrained
value, whereas $v$ is a candidate for the auxiliary value $V^\Psi$. The two are identified only after
verification; see \citet{WuEtAl2020,LiWangZhou2026}.}, let $v_t$, $v_x$, and $v_{xx}$ denote its time,
first-state, and second-state derivatives. The formal distribution-valued HJB equation is}
\begin{equation}
\left\{
\begin{aligned}
0={}&v_t(t,x)+\inf_{\pi\in\mathcal P_{\Psi,+}(t,x)}\Bigg\{
 v_x(t,x)\int b(t,x,u)\pi(du)\\
&+\frac12v_{xx}(t,x)\sum_{j=1}^m\int\sigma_j^2(t,x,u)\pi(du)
 +\int\ell(t,x,u)\pi(du)\Bigg\},\\
v(T,x)={}&q_Tx^2.
\end{aligned}
\right.
\label{eq:distributional_hjb}
\end{equation}

\begin{proposition}[Distribution-to-mean reduction]
\label{prop:distribution_to_mean}
Let $\pi\in\mathcal P_{\Psi,+}(t,x)$ and let $a:=\int_{\mathbb R^n}u\,\pi(du)\in\mathbb R_+^n$. Then
\begin{equation}
\int_{\mathbb R^n}uu^\top\pi(du)=aa^\top+x^2\Psi(t),
\label{eq:second_moment_reduction}
\end{equation}
\begin{equation}
\int_{\mathbb R^n}\ell(t,x,u)\pi(du)=a^\top R(t)a+2xS(t)^\top a+q(t)x^2+x^2\operatorname{tr}[R(t)\Psi(t)].
\label{eq:cost_reduction}
\end{equation}
{Together with Proposition~\ref{prop:exploratory_characteristics}, these identities imply that}
\eqref{eq:distributional_hjb} is reduced to the mean-valued HJB equation
\begin{equation}
\left\{
\begin{aligned}
0={}&v_t(t,x)+\inf_{a\in\mathbb R_+^n}\Bigg\{
 a^\top R(t)a+2xS(t)^\top a+q(t)x^2+x^2\operatorname{tr}[R(t)\Psi(t)]\\
&+\left[A(t)x+B(t)^\top a\right]v_x(t,x)\\
&+\frac12v_{xx}(t,x)\Big[\|xC(t)+D(t)^\top a\|^2
+x^2\operatorname{tr}\!\left[D(t)D(t)^\top\Psi(t)\right]\Big]\Bigg\},\\
v(T,x)={}&q_Tx^2.
\end{aligned}
\right.
\label{eq:mean_hjb}
\end{equation}
\end{proposition}

{The proof in \ref{app:main_proofs} also shows that every feasible mean
$a\in\mathbb R_+^n$ is attained by a probability measure with the prescribed covariance.}

{{With covariance fixed, the HJB integrand depends on the policy only through its
conditional mean.} Gaussianity is used later for computation, not for this reduction.}

\subsection{State separation and regional minimization}
\label{subsec:state_separation}

{
Following the regional quadratic representation in \citet{WuEtAl2020}, seek a candidate
for the auxiliary value of the form
\[
v(t,x)=
\begin{cases}
P_+^\Psi(t)x^2,&x>0,\\
0,&x=0,\\
P_-^\Psi(t)x^2,&x<0.
\end{cases}
\]
This is a candidate form; Theorem~\ref{thm:main} verifies it under the stated assumptions.
At a fixed time and in either open state region, write $p=P_\pm^\Psi(t)\ge0$.
Then $v_x=2px$ and $v_{xx}=2p$. Substituting these derivatives into
\eqref{eq:mean_hjb} gives the control-dependent terms
\[
a^\top\left[R(t)+pD(t)D(t)^\top\right]a
+2x\left\{S(t)+p[B(t)+D(t)C(t)]\right\}^\top a.
\]
For a generic $p\ge0$, collect these coefficients into the following matrices.
}
\begin{equation}
\widehat R(t,p):=R(t)+pD(t)D(t)^\top, \qquad
\widehat S(t,p):=S(t)+p\left[B(t)+D(t)C(t)\right].
\label{eq:regional_matrices}
\end{equation}
{For $p\ge0$, define the regional minimum values and feedback gains by}
\begin{equation}
{F_\pm(t,p):=\min_{K\in\mathbb R_+^n}}
\left\{K^\top\widehat R(t,p)K\pm2\widehat S(t,p)^\top K\right\},
\label{eq:regional_value}
\end{equation}
\begin{equation}
{K_\pm(t,p):=\arg\min_{K\in\mathbb R_+^n}}
\left\{K^\top\widehat R(t,p)K\pm2\widehat S(t,p)^\top K\right\}.
\label{eq:regional_gain}
\end{equation}
{Here the subscripts $+$ and $-$ refer to the positive and negative state regions. By
Assumption~\ref{ass:control_weight_positive_definite}, for every $p\ge0$,
\[
\widehat R(t,p)=R(t)+pD(t)D(t)^\top\succ0.
\]
Hence both regional objectives are strictly convex on $\mathbb R_+^n$, and the minimizers
$K_\pm(t,p)$ are uniquely defined.} Since $\mathbb R_+^n$ is a cone, scaling
$a=xK$ for $x>0$ and $a=(-x)K$ for $x<0$ gives
\begin{equation}
{\min}_{a\in\mathbb R_+^n}
\left\{a^\top\widehat R(t,p)a+2x\widehat S(t,p)^\top a\right\}
=
\begin{cases}
x^2F_+(t,p),&x>0,\\
0,&x=0,\\
x^2F_-(t,p),&x<0,
\end{cases}
\label{eq:state_separation_value}
\end{equation}
{The corresponding unique minimizer is}
\begin{equation}
a^*(t,x;p)=
\begin{cases}
xK_+(t,p),&x>0,\\
0,&x=0,\\
-xK_-(t,p),&x<0.
\end{cases}
\label{eq:state_separation_minimizer}
\end{equation}
{
For $x=0$, the linear term disappears and $\widehat R(t,p)\succ0$, so the minimum is attained uniquely
at $a=0$. The prescribed covariance is also zero. {By \eqref{eq:second_moment_reduction}, the optimal
randomized action therefore equals zero with probability one. We denote this probability law by $\delta_0$.}
} The calculations for all three
cases are recorded in \ref{app:main_proofs}.

{The constrained minimization thus separates into two finite-dimensional quadratic
programs, one for each open state region. Their optimal values determine the regional Riccati
equations.}

\subsection{Regional Riccati equations}
\label{subsec:generalized_riccati}

{
Use the regional quadratic candidate introduced in Section~\ref{subsec:state_separation}.
An overdot denotes differentiation with respect to time. Substitution into \eqref{eq:mean_hjb}
and the regional minimization give
}
\begin{equation}
\left\{
\begin{aligned}
0={}&\dot P_\pm^\Psi(t)+[2A(t)+\|C(t)\|^2]P_\pm^\Psi(t)+q(t)+\operatorname{tr}[R(t)\Psi(t)]\\
&+P_\pm^\Psi(t)\operatorname{tr}[D(t)D(t)^\top\Psi(t)]
+F_\pm\!\left(t,P_\pm^\Psi(t)\right),\\
P_\pm^\Psi(T)={}&q_T.
\end{aligned}
\right.
\label{eq:riccati_pm}
\end{equation}
{Define the regional gains on $[0,T]$ by}
\begin{equation}
K_\pm^\Psi(t):=K_\pm\!\left(t,P_\pm^\Psi(t)\right),
\qquad {t\in[0,T]}.
\label{eq:exploratory_regional_gains}
\end{equation}
{Assumption~\ref{ass:control_weight_positive_definite} ensures that these gains are uniquely defined for
every $t\in[0,T]$, including when $P_\pm^\Psi(t)=0$.} The derivation of \eqref{eq:riccati_pm} and the
uniqueness of the regional minimizers are given in \ref{app:main_proofs}.

{The two trace terms in \eqref{eq:riccati_pm} account for the prescribed covariance. They
can change $P_\pm^\Psi$ and, because $K_\pm(t,p)$ depends on $p$, can also change the optimal
conditional mean. The next theorem verifies the value candidate and identifies that mean.}

\subsection{Main theorem: value function and optimal conditional mean}
\label{subsec:main_theorem}

{The regional equations and state separation give the following characterization.}

\begin{theorem}[Value function and optimal conditional mean]
\label{thm:main}
Suppose {Assumptions~1--4} hold and $\Psi$ is deterministic, measurable, bounded on $[0,T]$, and positive
semidefinite. {Assume further that \eqref{eq:riccati_pm} admits nonnegative solutions
$P_+^\Psi,P_-^\Psi\in C^1([0,T])$.} Under the standard dynamic-programming and
comparison conditions for \eqref{eq:mean_hjb},
\begin{equation}
V^\Psi(t,x)=
\begin{cases}
P_+^\Psi(t)x^2,&x>0,\\
0,&x=0,\\
P_-^\Psi(t)x^2,&x<0,
\end{cases}
\label{eq:main_value}
\end{equation}
is the value function of the exploratory problem and the unique viscosity solution of
\eqref{eq:mean_hjb}. {The optimal conditional mean is}
\begin{equation}
a^{\Psi,*}(t,x)=
\begin{cases}
xK_+^\Psi(t),&x>0,\\
0,&x=0,\\
-xK_-^\Psi(t),&x<0.
\end{cases}
\label{eq:optimal_mean}
\end{equation}
\end{theorem}

\begin{proof}
The regional verification and the boundary viscosity argument are given in \ref{app:main_proofs}.
\end{proof}

{The theorem concerns the auxiliary problem; equality with the original optimal control
requires a further condition. {
Under the optimal feedback, both the mean and covariance vanish at $x=0$.
The action equals zero with probability one, so the state remains at zero.
} Away from this boundary, covariance can change the regional
gains. We next determine when it does not.}

\subsection{Control--cost compatibility and mean-policy invariance}
\label{subsec:alignment_consistency}

{To compare the two optimal feedbacks, set $\Psi\equiv0$. The original regional
coefficients $P_+^0$ and $P_-^0$ then satisfy}
\begin{equation}
\left\{
\begin{aligned}
0={}&\dot P_\pm^0(t)+[2A(t)+\|C(t)\|^2]P_\pm^0(t)+q(t)
+F_\pm\!\left(t,P_\pm^0(t)\right),\\
P_\pm^0(T)={}&q_T,
\end{aligned}
\right.
\label{eq:classical_riccati_pm}
\end{equation}
Let $u^*$ denote the optimal feedback control of the original problem. {Then}
\begin{equation}
K_\pm^0(t):=K_\pm\!\left(t,P_\pm^0(t)\right),
\qquad
u^*(t,x)=
\begin{cases}
xK_+^0(t),&x>0,\\
0,&x=0,\\
-xK_-^0(t),&x<0.
\end{cases}
\label{eq:classical_control_our_notation}
\end{equation}
In general, the exploration terms in \eqref{eq:riccati_pm} change $P_\pm^\Psi$, and because $K_\pm(t,p)$
may depend on $p$, the exploratory mean need not coincide with $u^*$.

{{
The compatibility condition below is motivated by the second-moment expansion of the
state equation.} Here $x$ denotes the scalar state without a particular financial interpretation.}

Over a short interval of length $\Delta t$, a direct expansion gives
\begin{align*}
\mathbb E_{t,x}[X(t+\Delta t)^2]-x^2
={}&\Big\{[2A(t)+\|C(t)\|^2]x^2\\
&\quad+u^\top D(t)D(t)^\top u\\
&\quad+2x[B(t)+D(t)C(t)]^\top u\Big\}\Delta t
+o(\Delta t).
\end{align*}
The first term is independent of the control. The control-dependent contribution is therefore
\[
\left\{
u^\top D(t)D(t)^\top u
+2x[B(t)+D(t)C(t)]^\top u
\right\}\Delta t.
\]
{The short-time conditional second-moment expansion shows that its control-dependent part is
determined by the quadratic coefficient $D(t)D(t)^\top$ and the state--control cross coefficient
$B(t)+D(t)C(t)$. We therefore require the running control cost to use these same two coefficients with
a common positive weight $\lambda(t)$. This requirement leads to the following compatibility condition.}

\begin{assumption}[Control--cost compatibility]
\label{ass:alignment}
There exists a deterministic measurable function {$\lambda:[0,T]\to(0,\infty)$} such that
\begin{equation}
R(t)=\lambda(t)D(t)D(t)^\top,
\qquad
S(t)=\lambda(t)[B(t)+D(t)C(t)],
\qquad t\in[0,T].
\label{eq:alignment_condition}
\end{equation}
\end{assumption}

{Assumption~\ref{ass:alignment} therefore identifies a structurally compatible subclass of the general LQ
formulation: the control-dependent running cost uses the same quadratic and cross coefficients as the
control-dependent part of the second-moment expansion. {The state-only penalty $q(t)x^2$
does not participate in this coefficient matching, but remains subject to the cost convexity in
Assumption~2.} The
assumption is used only to establish mean-policy invariance and is not required for the general
exploratory characterization obtained above.}

\begin{theorem}[Exact mean-policy consistency under control--cost compatibility]
\label{thm:alignment_consistency}
Suppose the assumptions of Theorem~\ref{thm:main} hold, {the original regional equations
\eqref{eq:classical_riccati_pm} admit nonnegative solutions}, and Assumption~\ref{ass:alignment} holds.
{Then, for every $t\in[0,T]$ and $p\ge0$,}
\begin{equation}
K_\pm(t,p)=\arg\min_{K\in\mathbb R_+^n}
\left\{K^\top D(t)D(t)^\top K
\pm2[B(t)+D(t)C(t)]^\top K\right\},
\label{eq:alignment_gain_independence}
\end{equation}
{Consequently, $K_\pm^\Psi(t)=K_\pm^0(t)$ for $t\in[0,T]$ and}
\begin{equation}
a^{\Psi,*}(t,x)=u^*(t,x),
\qquad {(t,x)\in[0,T]\times\mathbb R}.
\label{eq:alignment_mean_policy_identity}
\end{equation}
\end{theorem}

\begin{proof}
Assumption~\ref{ass:alignment} gives
\begin{equation}
\widehat R(t,p)=[p+\lambda(t)]D(t)D(t)^\top,
\qquad
\widehat S(t,p)=[p+\lambda(t)][B(t)+D(t)C(t)].
\label{eq:aligned_regional_matrices}
\end{equation}
{Because $p+\lambda(t)>0$ for every $p\ge0$, this positive common factor may be divided out without
changing the unique minimizer.} This proves \eqref{eq:alignment_gain_independence} and hence
\eqref{eq:alignment_mean_policy_identity}.
\end{proof}

{
Assumption~\ref{ass:alignment} is motivated by the preceding short-time expansion of the
state second moment. The expansion identifies $D(t)D(t)^\top$ and $B(t)+D(t)C(t)$ as the
coefficients of the control-dependent quadratic and cross terms. The assumption uses these
coefficients in the running control cost, weighted by $\lambda(t)$. The regional value coefficient
$p=P_\pm^\Psi(t)$ weights the current state squared in the optimal remaining cost, including
subsequent running costs and the terminal cost. In \eqref{eq:mean_hjb}, $p$ weights the same
control-dependent terms. Their weights add to the positive factor $p+\lambda(t)$ in
\eqref{eq:aligned_regional_matrices}, which leaves the minimizing gain unchanged. Gain independence
is therefore a consequence of this coefficient matching.}

\subsection{Infinite-horizon limit and mean-square stability}
\label{subsec:infinite_horizon}

{We now consider constant coefficients $A,B,C,D,R,S,q$, and
$\Psi\in\mathbb S_+^n$ on an infinite horizon. The randomized policy remains a general
{measurable feedback policy} with the moments specified in Section~\ref{subsec:randomized_policies}. We obtain the
{infinite-horizon problem} through zero-terminal-cost finite-horizon values.}

Extend the filtered probability space and Brownian motion in Section~\ref{subsec:notation} to $t\ge0$.
Define $\Pi_{\Psi,+}[0,\infty)$ to consist of {measurable randomized
feedback policies} $\pi_t(du\mid x)$ satisfying
\[
a^\pi(t,x)\in\mathbb R_+^n,\qquad
\Sigma^\pi(t,x)=x^2\Psi,
\]
and generating a state process on $[0,\infty)$ such that, for every $T>0$, the restriction $\pi|_{[0,T]}$
belongs to the finite-horizon admissible set $\Pi_{\Psi,+}[0,T]$ defined immediately after
\eqref{eq:exploratory_problem}. Thus admissibility is required on every finite horizon, but it does not
require the infinite-horizon cost to be finite.

With the time argument retained for a possibly time-dependent policy, set
\begin{equation}
\begin{aligned}
J_\infty^\Psi(x_0;\pi)
&:=\mathbb E\int_0^\infty
 \overline\ell^\pi(t,X^\pi(t))\,dt,\\
V_\infty^\Psi(x_0)
&:=\inf_{\pi\in\Pi_{\Psi,+}[0,\infty)}
 J_\infty^\Psi(x_0;\pi).
\end{aligned}
\label{eq:infinite_exploratory_problem}
\end{equation}
Here $\overline\ell^\pi(t,x)$ is the {running cost under $\pi$} defined in \eqref{eq:running_cost_under_policy}
with constant cost coefficients. The nonnegative integral is defined as the increasing limit of the
integrals over $[0,T]$ and may equal $+\infty$. {Finite-interval admissibility does not by
itself ensure $J_\infty^\Psi(x_0;\pi)<\infty$; we establish finite cost for the optimal candidate
below.}

\begin{assumption}[Infinite-horizon cost and diffusion conditions]
\label{ass:infinite_horizon_cost}
In addition to the constant-coefficient setting above,
\begin{equation}
Q=\begin{pmatrix}R&S\\S^\top&q\end{pmatrix}\succ0,
\qquad DD^\top\succ0.
\label{eq:infinite_cost_positive_definite}
\end{equation}
\end{assumption}

{These conditions strengthen the finite-horizon cost assumption, agree with the
diffusion condition in Assumption~4, and {leave the nonnegative gain set nonempty}.}

For $T>0$, consider the auxiliary finite-horizon problem
\begin{equation}
\begin{aligned}
J_T^\Psi(x_0;\pi)
&:=\mathbb E\int_0^T
 \overline\ell^\pi(t,X^\pi(t))\,dt,\\
V_T^\Psi(x_0)
&:=\inf_{\pi\in\Pi_{\Psi,+}[0,T]}J_T^\Psi(x_0;\pi),
\end{aligned}
\label{eq:infinite_auxiliary_finite_problem}
\end{equation}
subject to \eqref{eq:exploratory_state}, the mean constraint, and the prescribed covariance. This problem
has no terminal cost: $q_T=0$, so its continuation value at time $T$ is zero for every terminal state.
Consequently, both regional value coefficients satisfy
\[
P_+^\Psi(T;T)=P_-^\Psi(T;T)=0.
\]
{The zero terminal value does not impose $X^\pi(T)=0$.}

{For every $T>0$, the finite-horizon hypotheses of Theorem~\ref{thm:main} are additionally assumed; they
do not follow from \eqref{eq:infinite_cost_positive_definite}.}

Write $P_\pm^\Psi(t;T)$ for the regional coefficients associated with this auxiliary problem.
Theorem~\ref{thm:main} and \eqref{eq:riccati_pm}, applied with $q_T=0$, give
\begin{equation}
V_T^\Psi(x_0)
=x_0^2\left[
P_+^\Psi(0;T)\mathbf 1_{\{x_0\ge0\}}
+P_-^\Psi(0;T)\mathbf 1_{\{x_0<0\}}\right].
\label{eq:infinite_auxiliary_value}
\end{equation}
To write the coefficient equations explicitly, suppress the time argument in the already defined regional
programs \eqref{eq:regional_value}--\eqref{eq:regional_gain}. For $p\ge0$, write
\begin{equation}
\begin{aligned}
F_\pm(p)
&:=\min_{K\in\mathbb R_+^n}
 \left\{K^\top(R+pDD^\top)K
 \pm2[S+p(B+DC)]^\top K\right\},\\
K_\pm(p)
&:=\arg\min_{K\in\mathbb R_+^n}
 \left\{K^\top(R+pDD^\top)K
 \pm2[S+p(B+DC)]^\top K\right\}.
\end{aligned}
\label{eq:infinite_regional_qp}
\end{equation}
Then \eqref{eq:riccati_pm} becomes
\begin{equation}
\left\{
\begin{aligned}
0={}&\dot P_\pm^\Psi(t;T)
 +[2A+\|C\|^2+\operatorname{tr}(DD^\top\Psi)]P_\pm^\Psi(t;T)\\
&+q+\operatorname{tr}(R\Psi)
 +F_\pm(P_\pm^\Psi(t;T)),\qquad 0\le t\le T,\\
P_\pm^\Psi(T;T)={}&0.
\end{aligned}
\right.
\label{eq:infinite_finite_riccati}
\end{equation}
{The prescribed covariance produces the two trace terms through
Propositions~\ref{prop:exploratory_characteristics} and
\ref{prop:distribution_to_mean}. The regional minimizer $K_\pm(p)$ is unique because
$R+pDD^\top\succ0$.}

\begin{lemma}[Time translation of the auxiliary coefficients]
\label{lem:infinite_time_translation}
For $T>0$ and $\tau>0$, suppose the finite-horizon value characterization in Theorem~\ref{thm:main}
applies to the auxiliary problems with horizons $T$ and $T+\tau$. Then
\begin{equation}
P_\pm^\Psi(t;T)=P_\pm^\Psi(t+\tau;T+\tau),
\qquad 0\le t\le T.
\label{eq:infinite_time_translation}
\end{equation}
\end{lemma}

\begin{proof}
Write $V^\Psi(t,x;T)$ for the value of the zero-terminal-cost auxiliary problem on $[t,T]$. For a policy
$\pi$ on $[t,T]$, define its time shift on $[t+\tau,T+\tau]$ by
\[
\pi_{s+\tau}^{\tau}(du\mid x)=\pi_s(du\mid x),
\qquad t\le s\le T.
\]
The mean constraint and the covariance $x^2\Psi$ are unchanged, the Brownian increment processes have the
same law, and both problems have zero terminal cost. Hence the corresponding costs,
and therefore the value functions, are equal:
\[
V^\Psi(t,x;T)=V^\Psi(t+\tau,x;T+\tau).
\]
{Applying Theorem~\ref{thm:main} at $x=1$ and $x=-1$ gives
\eqref{eq:infinite_time_translation}.}
\end{proof}

\begin{theorem}[Infinite-horizon exploratory control]
\label{thm:infinite_horizon}
Suppose Assumptions~1--4, Assumption~\ref{ass:infinite_horizon_cost}, and the finite-horizon verification
conditions stated above hold in the constant-coefficient setting, and \eqref{eq:infinite_finite_riccati}
has nonnegative $C^1$ solutions for every finite horizon. Assume there is $M>0$ such that
\begin{equation}
P_+^\Psi(0;T)\le M,\qquad P_-^\Psi(0;T)\le M,
\qquad T>0.
\label{eq:infinite_uniform_bound}
\end{equation}
Then the finite limits
\begin{equation}
P_\pm^{\Psi,*}=\lim_{T\to\infty}P_\pm^\Psi(0;T)
\label{eq:infinite_coefficient_limits}
\end{equation}
exist. Moreover, if $P_+^\Psi(0;T)$ and $P_-^\Psi(0;T)$ are uniformly continuous with respect to
$T$, then $P_\pm^{\Psi,*}>0$ and they satisfy
\begin{equation}
\begin{aligned}
0={}&[2A+\|C\|^2+\operatorname{tr}(DD^\top\Psi)]P_\pm^{\Psi,*}\\
&+q+\operatorname{tr}(R\Psi)+F_\pm(P_\pm^{\Psi,*}).
\end{aligned}
\label{eq:infinite_algebraic_riccati}
\end{equation}
With the uniquely defined gains
\begin{equation}
K_\pm^{\Psi,*}=K_\pm(P_\pm^{\Psi,*}),
\label{eq:infinite_stationary_gains}
\end{equation}
set
\begin{equation}
a_\infty^{\Psi,*}(x)=
\begin{cases}
xK_+^{\Psi,*},&x>0,\\
0,&x=0,\\
-xK_-^{\Psi,*},&x<0.
\end{cases}
\label{eq:infinite_stationary_mean}
\end{equation}
{There exists an admissible randomized {feedback policy} $\pi(du\mid x)$ with this conditional
mean and covariance $x^2\Psi$, and every admissible randomized {feedback policy} with these moments is
optimal for \eqref{eq:infinite_exploratory_problem}.}
The value is
\begin{equation}
V_\infty^\Psi(x_0)
=x_0^2\left[
P_+^{\Psi,*}\mathbf 1_{\{x_0\ge0\}}
+P_-^{\Psi,*}\mathbf 1_{\{x_0<0\}}\right].
\label{eq:infinite_optimal_value}
\end{equation}
{Under each such policy, the closed-loop system is mean-square asymptotically stable:}
\begin{equation}
\lim_{t\to\infty}\mathbb E[|X^{\Psi,*}(t)|^2]=0.
\label{eq:infinite_mean_square_stability}
\end{equation}
\end{theorem}

\begin{proof}
For $T,\tau>0$, restriction of a policy in $\Pi_{\Psi,+}[0,T+\tau]$ to $[0,T]$ remains admissible. Since
$Q\succ0$, the {running cost under $\pi$} is nonnegative, and therefore $V_{T+\tau}^\Psi(x_0)\ge
V_T^\Psi(x_0)$. Evaluating \eqref{eq:infinite_auxiliary_value} at $x_0=1$ and $x_0=-1$ shows that
$P_+^\Psi(0;T)$ and $P_-^\Psi(0;T)$ are nondecreasing in $T$.
{The bound \eqref{eq:infinite_uniform_bound} then gives the limits in
\eqref{eq:infinite_coefficient_limits}. }

For the remaining conclusions, assume the stated uniform continuity and set
$g_\pm^\Psi(s)=P_\pm^\Psi(0;s)$ for $s>0$ and $g_\pm^\Psi(0)=0$. Here $s$ is the optimization horizon
length. Lemma~\ref{lem:infinite_time_translation} gives $P_\pm^\Psi(t;T)=g_\pm^\Psi(T-t)$; substituting
this identity into \eqref{eq:infinite_finite_riccati} gives
\begin{equation}
\left\{
\begin{aligned}
\dot g_\pm^\Psi(s)
={}&[2A+\|C\|^2+\operatorname{tr}(DD^\top\Psi)]g_\pm^\Psi(s)\\
&+q+\operatorname{tr}(R\Psi)+F_\pm(g_\pm^\Psi(s)),
\qquad s\ge0,\\
g_\pm^\Psi(0)={}&0.
\end{aligned}
\right.
\label{eq:infinite_forward_riccati}
\end{equation}
{The prescribed covariance enters this forward system through the two trace terms.}
Since $g_\pm^\Psi(s)\to P_\pm^{\Psi,*}$, $\int_0^s\dot
g_\pm^\Psi(v)\,dv=g_\pm^\Psi(s)-g_\pm^\Psi(0)$ has a finite limit. The parametric convex-program
continuity argument in \citet[footnote~12]{WuEtAl2020} applies to \eqref{eq:infinite_regional_qp}, which
has the same strongly convex quadratic structure; as in that footnote, it supplies the uniform continuity
of $\dot g_\pm^\Psi$ and the continuity of the optimized term. Barbalat's lemma, cited there from
\citet{Khalil1992}, gives $\dot g_\pm^\Psi(s)\to0$, and passage to the limit in
\eqref{eq:infinite_forward_riccati} yields \eqref{eq:infinite_algebraic_riccati}.

It remains to record the changes caused by prescribed covariance and to verify the randomized-policy
conclusion. For a constant regional gain $K\in\mathbb R_+^n$,
Propositions~\ref{prop:exploratory_characteristics} and \ref{prop:distribution_to_mean} give the
second-moment growth rate and the coefficient of $x^2$ in the {running cost under $\pi$} as
\begin{equation}
\begin{aligned}
N_\pm^\Psi(K)
&:=2A+\|C\|^2+\operatorname{tr}(DD^\top\Psi)\\
&\quad+K^\top DD^\top K\pm2(B+DC)^\top K,\\
L_\pm^\Psi(K)
&:=q+\operatorname{tr}(R\Psi)+K^\top RK\pm2S^\top K.
\end{aligned}
\label{eq:infinite_rates}
\end{equation}
{The first trace term in \eqref{eq:infinite_rates} reflects exploratory diffusion;
the second reflects exploratory running cost. At $K=K_\pm^{\Psi,*}$,
\eqref{eq:infinite_algebraic_riccati} gives}
\begin{equation}
P_\pm^{\Psi,*}N_\pm^\Psi(K_\pm^{\Psi,*})
=-L_\pm^\Psi(K_\pm^{\Psi,*}).
\label{eq:infinite_rate_identity}
\end{equation}
Indeed, $Q\succ0$ and $\Psi\succeq0$ imply
\[
L_\pm^\Psi(K)=
\begin{pmatrix}K\\ \pm1\end{pmatrix}^{\!\top}
Q\begin{pmatrix}K\\ \pm1\end{pmatrix}
+\operatorname{tr}(R\Psi)>0.
\]
If $P_\pm^{\Psi,*}=0$, \eqref{eq:infinite_algebraic_riccati} would instead give $0=\min_{K\in\mathbb
R_+^n}L_\pm^\Psi(K)>0$. Hence $P_\pm^{\Psi,*}>0$, and \eqref{eq:infinite_rate_identity} yields
\begin{equation}
N_\pm^\Psi(K_\pm^{\Psi,*})
=-\frac{L_\pm^\Psi(K_\pm^{\Psi,*})}{P_\pm^{\Psi,*}}<0.
\label{eq:infinite_negative_rates}
\end{equation}

{These quantities depend only on the first two conditional moments and therefore apply
to every {admissible randomized {feedback policy}} with mean \eqref{eq:infinite_stationary_mean} and covariance
$x^2\Psi$. For $x_0\ne0$, set $\varepsilon=+$ when $x_0>0$ and $\varepsilon=-$ when $x_0<0$.
The scalar homogeneous closed loop preserves the initial half-line. It\^o's formula then gives}
\begin{equation}
\frac{d}{dt}\mathbb E[|X^{\Psi,*}(t)|^2]
=N_\varepsilon^\Psi(K_\varepsilon^{\Psi,*})
 \mathbb E[|X^{\Psi,*}(t)|^2],
\qquad
\mathbb E[|X^{\Psi,*}(0)|^2]=x_0^2.
\label{eq:infinite_second_moment_ode}
\end{equation}
Solving this scalar differential equation yields
\begin{equation}
\mathbb E[|X^{\Psi,*}(t)|^2]
=\begin{cases}
x_0^2\exp(N_+^\Psi(K_+^{\Psi,*})t),&x_0>0,\\
0,&x_0=0,\\
x_0^2\exp(N_-^\Psi(K_-^{\Psi,*})t),&x_0<0.
\end{cases}
\label{eq:infinite_second_moment_solution}
\end{equation}
For $x_0=0$, the state remains zero. The negative rates in \eqref{eq:infinite_negative_rates} therefore
prove \eqref{eq:infinite_mean_square_stability}; the exponential form follows from the displayed
constant-coefficient second-moment equation and does not use Gaussianity.
\end{proof}

{The uniform bound \eqref{eq:infinite_uniform_bound} is an additional assumption; neither
monotonicity nor \eqref{eq:infinite_cost_positive_definite} implies it.\footnote{{For
example, let $n=m=1$, $A=1$, $B=C=S=0$, $D=R=q=1$, and $\Psi=\psi\ge0$. Then
$Q=I_2\succ0$, and the regional minimization reduces, for either sign, to
$\min_{K\ge0}(1+p)K^2$, whose unique minimizer is $K=0$. Equation
\eqref{eq:infinite_forward_riccati} becomes
$\dot g(s)=(2+\psi)g(s)+(1+\psi)$ with $g(0)=0$, so
$g(s)=\frac{1+\psi}{2+\psi}(e^{(2+\psi)s}-1)\to\infty$. Thus the uniform bound fails, and the
associated algebraic equation $0=(2+\psi)P+(1+\psi)$ has no nonnegative solution.}}}

{The finite-horizon compatibility condition also yields {an infinite-horizon mean-policy
identity}.}

\begin{proposition}[Infinite-horizon mean-policy consistency]
\label{prop:infinite_alignment}
Suppose the assumptions of Theorem~\ref{thm:infinite_horizon} hold both for the exploratory problem with
covariance $\Psi$ and for its original-control specialization $\Psi=0$. If
the time-invariant compatibility condition
\begin{equation}
R=\lambda DD^\top,
\qquad
S=\lambda(B+DC),
\qquad \lambda>0,
\label{eq:infinite_alignment}
\end{equation}
holds, then
\begin{equation}
K_\pm^{\Psi,*}=K_\pm^{0,*}
=\arg\min_{K\in\mathbb R_+^n}
\left\{K^\top DD^\top K
\pm2(B+DC)^\top K\right\}.
\label{eq:infinite_gain_invariance}
\end{equation}
Consequently, {the optimal conditional mean of the exploratory problem coincides with the
optimal feedback of the {original infinite-horizon problem with nonnegative controls}}.
\end{proposition}

\begin{proof}
Equation \eqref{eq:infinite_alignment} gives the time-invariant version of
\eqref{eq:aligned_regional_matrices}. For every $p\ge0$, the objective in \eqref{eq:infinite_regional_qp}
is the positive factor $p+\lambda$ times the objective displayed in \eqref{eq:infinite_gain_invariance}.
Its unique minimizer is therefore independent of both $p$ and $\Psi$, by the same argument as in
Theorem~\ref{thm:alignment_consistency}. Applying this identity at
$p=P_\pm^{\Psi,*}$ and $p=P_\pm^{0,*}$ proves the result.
\end{proof}

\section{Solving the Constrained Stochastic LQ Problem with Unknown Dynamics}
\label{sec:unknown_dynamics}

{{Unknown dynamics prevent direct computation of the feedback gains
$K_+^\Psi(t)$ and $K_-^\Psi(t)$.} We therefore
learn the exploratory mean from sampled trajectories. Under
Assumption~\ref{ass:alignment}, its optimal target is also the original constrained feedback by
Theorem~\ref{thm:alignment_consistency}; otherwise the target is the exploratory optimum in
Theorem~\ref{thm:main}.}

\subsection{Policy parameterization and evaluation}
\label{subsec:rl_policy_evaluation}

Fix $\Psi_0\in\mathbb S_{++}^n$ and set
\begin{equation}
\Psi(t)\equiv\Psi_0.
\label{eq:rl_fixed_covariance}
\end{equation}
{The fixed matrix $\Psi_0$ sets the action covariance, while the conditional mean remains
to be learned. Although a rectified linear unit output could enforce a nonnegative mean, the regional
structure in Theorem~\ref{thm:main} suggests the parameterization}
\begin{equation}
a^\phi(t,x)=
\begin{cases}
x\phi_+(t),&x>0,\\
0,&x=0,\\
-x\phi_-(t),&x<0,
\end{cases}
\qquad
\phi_+(t),\phi_-(t)\in\mathbb R_+^n,
\label{eq:rl_param_mean}
\end{equation}
and use the Gaussian computational policy
\begin{equation}
\pi_t^\phi(du\mid x)=
\begin{cases}
\mathcal N\!\left(a^\phi(t,x),x^2\Psi_0\right)(du),&x\ne0,\\
\delta_0(du),&x=0.
\end{cases}
\label{eq:rl_gaussian_policy}
\end{equation}
{This Gaussian family is a computational representative of the admissible laws in
Section~\ref{subsec:randomized_policies}. {Its mean satisfies the nonnegative control constraint} and its
covariance is $x^2\Psi_0$.}

For a positive integer $N$, discretize $[0,T]$ by $t_i=iT/N$,
$i=0,\ldots,N$, with $\Delta t=T/N$, and let
\[
\phi_{+,i}:=\phi_+(t_i),
\qquad
\phi_{-,i}:=\phi_-(t_i).
\]
Thus the actor has $2N+2$ vector parameters. We use $\phi$ for both $(\phi_+,\phi_-)$ and the
sequence $\{(\phi_{+,i},\phi_{-,i})\}_{i=0}^N$.

For a fixed policy $\pi^\phi$, let $\mathbb E_{t,x}^{\pi^\phi}$ denote expectation for the trajectory
generated by $\pi^\phi$ and initialized by $X(t)=x$. Define
\begin{align}
V^\phi(t,x)
&:=\mathbb E_{t,x}^{\pi^\phi}
\left[\int_t^T\overline\ell^\phi(s,X(s))ds+q_TX(T)^2\right],\nonumber\\
\overline\ell^\phi(t,x)
&:=\int_{\mathbb R^n}\ell(t,x,u)\pi_t^\phi(du\mid x).
\label{eq:rl_fixed_policy_value}
\end{align}
On each open region, $v:=V^\phi$ satisfies the Feynman--Kac equation
\begin{equation}
\left\{
\begin{aligned}
0={}&v_t(t,x)+\left[A(t)x+B(t)^\top a^\phi(t,x)\right]v_x(t,x)\\
&+\frac12v_{xx}(t,x)
\left\|xC(t)+D(t)^\top a^\phi(t,x)\right\|^2\\
&+\frac12v_{xx}(t,x)x^2
\operatorname{tr}\!\left(D(t)D(t)^\top\Psi_0\right)
+\overline\ell^\phi(t,x),\\
v(T,x)={}&q_Tx^2.
\end{aligned}
\right.
\label{eq:rl_fixed_policy_pde}
\end{equation}
{A piecewise-quadratic fixed-policy value takes the form}
\begin{equation}
V^\phi(t,x)=
\begin{cases}
\theta_+^\phi(t)x^2,&x>0,\\
0,&x=0,\\
\theta_-^\phi(t)x^2,&x<0,
\end{cases}
\label{eq:rl_candidate_fixed_policy_value}
\end{equation}
where the theoretical regional coefficients $\theta_+^\phi,\theta_-^\phi:[0,T]\to\mathbb R_+$ satisfy
\begin{align}
0={}&\dot\theta_+^\phi
+\theta_+^\phi\left[
2A+2B^\top\phi_+
+\|C+D^\top\phi_+\|^2
+\operatorname{tr}(DD^\top\Psi_0)
\right]\nonumber\\
&+\phi_+^\top R\phi_+
+2S^\top\phi_+
+q+\operatorname{tr}(R\Psi_0),
\qquad
\theta_+^\phi(T)=q_T,
\label{eq:rl_positive_theta_equation}\\
0={}&\dot\theta_-^\phi
+\theta_-^\phi\left[
2A-2B^\top\phi_-
+\|C-D^\top\phi_-\|^2
+\operatorname{tr}(DD^\top\Psi_0)
\right]\nonumber\\
&+\phi_-^\top R\phi_-
-2S^\top\phi_-
+q+\operatorname{tr}(R\Psi_0),
\qquad
\theta_-^\phi(T)=q_T.
\label{eq:rl_negative_theta_equation}
\end{align}
All time arguments in \eqref{eq:rl_positive_theta_equation}--\eqref{eq:rl_negative_theta_equation} are
suppressed for readability; Appendix~\ref{app:rl_derivations} gives the substitution.

Because the dynamics are unknown, parameterize these coefficients by functions $\theta_+$ and
$\theta_-$ represented on the time grid by
\[
\theta_{+,i}:=\theta_+(t_i),
\qquad
\theta_{-,i}:=\theta_-(t_i),
\qquad
\theta_{+,i},\theta_{-,i}\in\mathbb R_+,
\qquad i=0,\ldots,N,
\]
and use $\theta$ for both the pair $(\theta_+,\theta_-)$ and the sequence
$\{(\theta_{+,i},\theta_{-,i})\}_{i=0}^N$. The critic is
\begin{equation}
V^\theta(t_i,x)=
\begin{cases}
\theta_{+,i}x^2,&x>0,\\
0,&x=0,\\
\theta_{-,i}x^2,&x<0.
\end{cases}
\label{eq:rl_param_value}
\end{equation}

Let $(X,U)$ denote a state--action trajectory generated by $\pi^\phi$, where $U(t)$ is sampled from
$\pi_t^\phi(\cdot\mid X(t))$. Write $X_i:=X(t_i)$ and $U_i:=U(t_i)$, and define
\begin{equation}
G_i:=\sum_{k=i}^{N-1}\ell(t_k,X_k,U_k)\Delta t+q_TX_N^2,
\qquad i=0,\ldots,N-1,
\label{eq:rl_realized_return}
\end{equation}
Thus $G_i$ is the realized cost-to-go from $t_i$ along one sampled trajectory. Its conditional
expectation under the fixed policy is approximated by the critic
$V^\theta(t_i,X_i)$. The corresponding critic loss is
\begin{equation}
L(\theta):=
\mathbb E\left[
\sum_{i=0}^{N-1}
\left(V^\theta(t_i,X_i)-G_i\right)^2\Delta t
\right].
\label{eq:rl_martingale_loss}
\end{equation}
{This is the martingale loss criterion of \citet[Section~3.1]{JiaZhou2022PE}, also used
for finite-horizon policy evaluation by \citet[Section~4.1]{LiWangZhou2026}. Each $G_i$ requires a
trajectory through $T$, so we approximate the loss with batches of complete episodes. The regional
gradient contributions from one trajectory are}
\begin{equation}
\begin{aligned}
g_{\theta,+,i}
&:=
\begin{cases}
2\Delta t\left[V^\theta(t_i,X_i)-G_i\right]X_i^2,&X_i>0,\\
0,&X_i\le0,
\end{cases}\\[2mm]
g_{\theta,-,i}
&:=
\begin{cases}
0,&X_i\ge0,\\
2\Delta t\left[V^\theta(t_i,X_i)-G_i\right]X_i^2,&X_i<0.
\end{cases}
\end{aligned}
\label{eq:rl_critic_gradient}
\end{equation}
For a mini-batch of $M$ trajectories, with superscript $(j)$ indexing the trajectory, use
\begin{equation}
\widehat g_{\theta,+,i}:=\frac1M\sum_{j=1}^M g_{\theta,+,i}^{(j)},
\qquad
\widehat g_{\theta,-,i}:=\frac1M\sum_{j=1}^M g_{\theta,-,i}^{(j)}.
\label{eq:rl_critic_sample_gradient}
\end{equation}
{For a vector argument, $\max\{0,\cdot\}$ is understood componentwise. With
critic learning rate $\alpha_\theta>0$, update}
\begin{equation}
\begin{aligned}
\theta_{+,i}&\leftarrow
{\max\!\left\{0,\theta_{+,i}-\alpha_\theta\widehat g_{\theta,+,i}\right\}},\\
\theta_{-,i}&\leftarrow
{\max\!\left\{0,\theta_{-,i}-\alpha_\theta\widehat g_{\theta,-,i}\right\}}.
\end{aligned}
\label{eq:rl_critic_update}
\end{equation}
The terminal parameters are fixed at
\begin{equation}
\phi_{+,N}=\phi_{-,N}=0,
\qquad
\theta_{+,N}=\theta_{-,N}=q_T.
\label{eq:rl_terminal_parameters}
\end{equation}

{The critic approximates the value of the current policy. We use it next to construct a
trajectory-based actor update.}

\subsection{Policy gradient}
\label{subsec:rl_policy_gradient}

{We use a continuous-time trajectory representation of the policy gradient. On either
open state region, $\pi_t^\phi(u\mid x)$ denotes the Gaussian density in
\eqref{eq:rl_gaussian_policy}. Assume that $\pi^\phi$ is admissible and smooth in $\phi$, its score is square integrable and continuous in the state--time variables, the policy-gradient integrand
is integrable, and the stochastic integrals below are true martingales.}

\begin{lemma}[Policy-gradient representation]
\label{lem:rl_policy_gradient}
On either open state region,
\begin{align}
\nabla_\phi V^\phi(t,x)
={}&\mathbb E_{t,x}\bigg[
\int_t^T
\nabla_\phi\log\pi_s^\phi(U(s)\mid X(s))\nonumber\\
&\qquad\qquad\times
\left\{
dV^\phi(s,X(s))
+\ell(s,X(s),U(s))ds
\right\}
\bigg].
\label{eq:rl_policy_gradient_representation}
\end{align}
\end{lemma}
\begin{proof}
The result follows directly from Theorem~5 of \citet{JiaZhou2022PG}.
\end{proof}

{Lemma~\ref{lem:rl_policy_gradient} expresses the fixed-policy value gradient through
observable trajectories rather than the unknown coefficients. Because the objective includes a running
LQ cost, the increment contains both the value change and that cost. {On the grid
$\{t_i\}_{i=0}^N$, replace $V^\phi$ in the policy-gradient representation by the critic $V^\theta$ and
approximate the running-cost integral over $[t_i,t_{i+1}]$ by its left-endpoint value. This gives the
observable increment}}
\begin{equation}
\Delta M_i^\theta
:=
V^\theta(t_{i+1},X_{i+1})
-V^\theta(t_i,X_i)
+\ell(t_i,X_i,U_i)\Delta t.
\label{eq:rl_martingale_increment}
\end{equation}
{For \eqref{eq:rl_gaussian_policy}, the score for each regional mean parameter is the
covariance-weighted deviation of the normalized action from its mean. Multiplying the score by
$\Delta M_i^\theta$ gives the one-trajectory regional contributions{\footnote{In the actor step, $\theta$ denotes the updated critic $\theta^{(k+1)}$, where $k$
indexes training rounds. Both value terms in \eqref{eq:rl_martingale_increment} use this critic,
with their respective time-grid coefficients. The score uses the actor that generated the retained
batch.}}}
\begin{equation}
\begin{aligned}
g_{\phi,+,i}
&:=
\begin{cases}
\displaystyle
\Psi_0^{-1}\left(\frac{U_i}{X_i}-\phi_{+,i}\right)\Delta M_i^\theta,
&X_i>0,\\
0,&X_i\le0,
\end{cases}\\[2mm]
g_{\phi,-,i}
&:=
\begin{cases}
0,&X_i\ge0,\\
\displaystyle
\Psi_0^{-1}\left(\frac{U_i}{-X_i}-\phi_{-,i}\right)\Delta M_i^\theta,
&X_i<0.
\end{cases}
\end{aligned}
\label{eq:rl_actor_gradient_contributions}
\end{equation}
For a mini-batch of $M$ trajectories, set
\begin{equation}
\widehat g_{\phi,+,i}:=\frac1M\sum_{j=1}^M g_{\phi,+,i}^{(j)},
\qquad
\widehat g_{\phi,-,i}:=\frac1M\sum_{j=1}^M g_{\phi,-,i}^{(j)}.
\label{eq:rl_actor_gradient}
\end{equation}
{With actor learning rate $\alpha_\phi>0$, use the nonnegativity-preserving update}
\begin{equation}
\begin{aligned}
\phi_{+,i}&\leftarrow
{\max\!\left\{0,\phi_{+,i}-\alpha_\phi\widehat g_{\phi,+,i}\right\}},\\
\phi_{-,i}&\leftarrow
{\max\!\left\{0,\phi_{-,i}-\alpha_\phi\widehat g_{\phi,-,i}\right\}}.
\end{aligned}
\label{eq:rl_actor_update}
\end{equation}

\begin{algorithmblock}{{Finite-horizon model-free actor--critic scheme}}
\label{alg:rl_model_free}

\textbf{Input.}
Choose $x_0\ne0$, the grid $\{t_i\}_{i=0}^N$, $\Psi_0$, initial parameters $\phi$ and $\theta$, learning
rates $\alpha_\phi,\alpha_\theta>0$, batch size $M$, and iteration numbers
$n_{\mathrm{warm}}$, $n_\theta$, and $n_\phi$.

\medskip
\textbf{Critic warm-up.}
For $k=1,\ldots,n_{\mathrm{warm}}$, sample $M$ trajectories under $\pi^\phi$, compute $G_i$, and update
$\theta$ using \eqref{eq:rl_critic_sample_gradient}--\eqref{eq:rl_critic_update}.

\medskip
\textbf{Actor--critic iterations.}
For $k=1,\ldots,n_\phi$:
\begin{enumerate}
\item {
Hold the actor fixed, perform $n_\theta$ critic updates, and retain the last sampled batch.}
\item {
On the retained batch, evaluate \eqref{eq:rl_martingale_increment} with the updated critic.}
\item {
Estimate \eqref{eq:rl_actor_gradient} with the actor that generated this batch and
update $\phi$ by \eqref{eq:rl_actor_update}.
}
\end{enumerate}

\textbf{Output.}
Return the updated policy $\pi^\phi$.
\end{algorithmblock}

{The algorithm returns an exploratory policy estimate and its conditional mean. Under the
state-scaled parameterization, the drift and diffusion vanish at $x=0$, and a trajectory from a nonzero
state stays in its initial region. A data set from one initial region therefore updates only its actor
and critic. Under Assumption~\ref{ass:alignment}, the target mean gain is also the corresponding gain
of the original constrained problem.}

\subsection{{
Infinite-horizon actor--critic algorithm
}}
\label{subsec:rl_infinite_horizon}

{The undiscounted infinite-horizon objective in
\eqref{eq:infinite_exploratory_problem} requires a different policy-evaluation signal. A complete
finite-horizon episode supplies the return $G_i$, whereas an infinite-horizon trajectory has no
terminal return. We therefore evaluate {the infinite-horizon policy} through martingale orthogonality
conditions on observed increments, following the policy-evaluation characterization of
\citet[Section~3.2]{JiaZhou2022PE}.}

{At each iteration, independent paths share a regional actor and critic, which are updated
in policy-evaluation--policy-improvement order. The actor uses a score-weighted one-step increment of
the type developed by \citet[Algorithm~3]{JiaZhou2022PG}; here the value is {an infinite-horizon total-cost
value} and the updates minimize cost.}

{The scheme uses separate positive- and negative-region parameters and projects both mean
gains onto $\mathbb R_+^n$. Each training round uses finite observation windows to evaluate and then
improve the shared regional policy, without imposing a finite-horizon terminal objective.}

Throughout this subsection the coefficients are constant, Assumption~\ref{ass:infinite_horizon_cost}
continues to hold, and $\Psi(t)\equiv\Psi_0\in\mathbb S_{++}^n$. Parameterize the {time-independent conditional
mean} and its Gaussian computational representative by
\begin{equation}
a^K(x)=
\begin{cases}
xK_+,&x>0,\\
0,&x=0,\\
-xK_-,&x<0,
\end{cases}
\qquad K_+,K_-\in\mathbb R_+^n,
\label{eq:rl_infinite_stationary_mean}
\end{equation}
and
\begin{equation}
\pi^K(du\mid x)=
\begin{cases}
\mathcal N(xK_+,x^2\Psi_0)(du),&x>0,\\
\delta_0(du),&x=0,\\
\mathcal N((-x)K_-,x^2\Psi_0)(du),&x<0.
\end{cases}
\label{eq:rl_infinite_gaussian_policy}
\end{equation}
The projection of $K_+$ and $K_-$ onto $\mathbb R_+^n$ below enforces the same conditional-mean
constraint as \eqref{eq:mean_noshorting}; individual Gaussian action
samples need not be componentwise nonnegative.

{For fixed $K=(K_+,K_-)$, use the regional rates $N_\pm^{\Psi_0}(K_\pm)$ and costs
$L_\pm^{\Psi_0}(K_\pm)$ from \eqref{eq:infinite_rates}. In the region
$\varepsilon\in\{+,-\}$, the second moment is
$\mathbb E_x^K[X(t)^2]=x^2\exp\{N_\varepsilon^{\Psi_0}(K_\varepsilon)t\}$. Thus, when}
\begin{equation}
N_\varepsilon^{\Psi_0}(K_\varepsilon)<0,
\label{eq:rl_infinite_fixed_policy_stability}
\end{equation}
the fixed-policy total cost and its boundary condition at infinity are
\begin{equation}
V_\infty^K(x)=\theta_\varepsilon^Kx^2,
\qquad
\theta_\varepsilon^K
=-\frac{L_\varepsilon^{\Psi_0}(K_\varepsilon)}
        {N_\varepsilon^{\Psi_0}(K_\varepsilon)}>0.
\label{eq:rl_infinite_fixed_policy_value}
\end{equation}
Moreover,
\begin{equation}
\lim_{T\to\infty}\mathbb E_x^K
   \bigl[V_\infty^K(X(T))\bigr]=0.
\label{eq:rl_infinite_transversality}
\end{equation}
{ In region
$\varepsilon$, its integrand is
$L_\varepsilon^{\Psi_0}(K_\varepsilon)\mathbb E_x^K[X(t)^2]$. If
$N_\varepsilon^{\Psi_0}(K_\varepsilon)\ge0$, the positive running-cost coefficient instead makes the
total cost infinite. The updates below use finite observed transitions without computing the rate
$N_\pm^{\Psi_0}$. At an intermediate iterate, their total-cost interpretation therefore requires
\eqref{eq:rl_infinite_fixed_policy_stability}; the algorithm does not guarantee stability at every
iterate. The optimal target is stabilizing by Theorem~\ref{thm:infinite_horizon}.}

{Use the regional quadratic critic}
\begin{equation}
V_\infty^\theta(x)=
\begin{cases}
\theta_+x^2,&x>0,\\
0,&x=0,\\
\theta_-x^2,&x<0,
\end{cases}
\qquad \theta_+,\theta_-\in\mathbb R_+.
\label{eq:rl_infinite_critic}
\end{equation}
For fixed actor and critic parameters, write
\begin{equation}
M_t^{\theta,K}
:=V_\infty^\theta(X(t))+\int_0^t\ell(X(s),U(s))ds.
\label{eq:rl_infinite_value_process}
\end{equation}
For an observed transition $(X_i,U_i,X_{i+1})$ over a time interval $\Delta t$, {use the following
one-step observed increment of the value-plus-accumulated-cost process:}
\begin{equation}
\Delta M_{\infty,i}^\theta
:=V_\infty^\theta(X_{i+1})-V_\infty^\theta(X_i)
  +\ell(X_i,U_i)\Delta t.
\label{eq:rl_infinite_martingale_increment}
\end{equation}
{The increment in \eqref{eq:rl_infinite_martingale_increment} has no terminal, discount,
or average-reward correction. If the gain satisfies
\eqref{eq:rl_infinite_fixed_policy_stability}, the fixed-policy calculation and It\^o's formula make
$M^{\theta^K,K}$ a martingale under the integrability conditions of
\citet{JiaZhou2022PE}.}

{Use the regional critic features as test functions in the martingale orthogonality
conditions:}
\begin{equation}
\xi_{+,i}=X_i^2\mathbf 1_{\{X_i>0\}},
\qquad
\xi_{-,i}=X_i^2\mathbf 1_{\{X_i<0\}},
\label{eq:rl_infinite_martingale_features}
\end{equation}
{which are the regional components of $\partial V_\infty^\theta(X_i)/\partial\theta$.
At round $k$, collect $M$ independent paths
$\{(X_i^{(m)},U_i^{(m)},X_{i+1}^{(m)})\}_{i=0}^{J-1}$ from the same nonzero regional initial state
$x_0^\varepsilon$. Each path uses $K_{\varepsilon,k}$ throughout the $J$ observed cells; superscript
$(m)$ indexes paths. The critic averages the pathwise sums of the orthogonality increments:}
{
\begin{equation}
\begin{aligned}
\widehat G_{\theta,\varepsilon,k}
&:=\frac1M\sum_{m=1}^M\sum_{i=0}^{J-1}
  \xi_{\varepsilon,i}^{(m)}
  \Delta M_{\infty,i}^{\theta_{\varepsilon,k},(m)},\\
\theta_{\varepsilon,k+1}
&:=\max\bigl\{0,\theta_{\varepsilon,k}
       +\alpha_{\theta,k}\widehat G_{\theta,\varepsilon,k}\bigr\}.
\end{aligned}
\label{eq:rl_infinite_critic_update}
\end{equation}
}

\begin{samepage}
{Following the regional score representation in \citet[equation~(26)]{LiWangZhou2026},
define the score-weighted actor contributions for an individual transition by{\footnote{In the actor step, $\Delta M_{\infty,i}^\theta$ is recomputed by
\eqref{eq:rl_infinite_martingale_increment} on the same batch with $\theta_{\varepsilon,k+1}$ after
\eqref{eq:rl_infinite_critic_update}. Both state values use the updated critic; the score uses
$K_{\varepsilon,k}$.}}}
{
\begin{equation}
\begin{aligned}
\Delta K_{+,i}^{\theta,K}
&:=
\begin{cases}
\displaystyle
\Psi_0^{-1}\left(\frac{U_i}{X_i}-K_+\right)
\Delta M_{\infty,i}^\theta,&X_i>0,\\
0,&X_i\le0,
\end{cases}\\[2mm]
\Delta K_{-,i}^{\theta,K}
&:=
\begin{cases}
0,&X_i\ge0,\\
\displaystyle
\Psi_0^{-1}\left(\frac{U_i}{-X_i}-K_-\right)
\Delta M_{\infty,i}^\theta,&X_i<0.
\end{cases}
\end{aligned}
\label{eq:rl_infinite_actor_increment}
\end{equation}
}
\end{samepage}
{
\begin{equation}
\begin{aligned}
\widehat G_{K,\varepsilon,k}
&:=\frac1M\sum_{m=1}^M\sum_{i=0}^{J-1}
  \left(\Delta K_{\varepsilon,i}^{
    \theta_{\varepsilon,k+1},K_{\varepsilon,k}}\right)^{(m)},\\
K_{\varepsilon,k+1}
&:=\max\bigl\{0,K_{\varepsilon,k}
       -\alpha_{K,k}\widehat G_{K,\varepsilon,k}\bigr\}.
\end{aligned}
\label{eq:rl_infinite_actor_update}
\end{equation}
}
{Equation~\eqref{eq:rl_infinite_actor_update} projects the resulting cost-descent step onto
the nonnegative cone.}

% {
% \begin{proposition}[Population consistency of the infinite-horizon updates]
% \label{prop:rl_infinite_population_consistency}
% Fix a region $\varepsilon\in\{+,-\}$ and a gain satisfying
% \eqref{eq:rl_infinite_fixed_policy_stability}. The critic moment has zero precisely
% when
% \begin{equation}
% \theta_\varepsilon N_\varepsilon^{\Psi_0}(K_\varepsilon)
% +L_\varepsilon^{\Psi_0}(K_\varepsilon)=0.
% \label{eq:rl_infinite_critic_population_equation}
% \end{equation}
% The corresponding actor moment is
% \begin{equation}
% \begin{aligned}
% \frac{1}{x^2dt}\mathbb E
%  \bigl[\mathcal S_{\pm,t}^K\,dM_t^{\theta,K}\mid X(t)=x\bigr]
% =2\bigl[&(R+\theta_\pm DD^\top)K_\pm\\
% &\ \pm\{S+\theta_\pm(B+DC)\}\bigr].
% \end{aligned}
% \label{eq:rl_infinite_actor_population_equation}
% \end{equation}
% Hence the critic zero and the projected actor zero coincide with the fixed-policy value equation and the
% first-order condition of the regional quadratic program
% \eqref{eq:infinite_regional_qp}. Their joint population fixed point solves
% \eqref{eq:infinite_algebraic_riccati}.
% \end{proposition}

% \begin{proof}
% The first identity follows from It\^o's formula applied to
% \eqref{eq:rl_infinite_value_process}; the second follows from the centered Gaussian
% moments of the score in \eqref{eq:rl_infinite_actor_increment}. The calculation is given in
% Appendix~\ref{app:rl_infinite_population}. The final statement follows by
% comparison with \eqref{eq:infinite_regional_qp}.
% \end{proof}
% }

{
\begin{algorithmblock}{{Batch actor--critic for the undiscounted infinite horizon}}
\label{alg:rl_infinite_actor_critic}

\textbf{Input.}
Choose $\Psi_0\in\mathbb S_{++}^n$, observation step $\Delta t$, path count $M$, observation-window
cell count $J$, regional initial states $x_0^+>0$ and $x_0^-<0$, initial parameters
$\theta_\pm\ge0$ and $K_\pm\in\mathbb R_+^n$, and learning-rate sequences
$\{\alpha_{\theta,k}\}$ and $\{\alpha_{K,k}\}$.

\medskip
\textbf{Learning.}
For $k=0,1,\ldots$ and each region $\varepsilon\in\{+,-\}$:
\begin{enumerate}
\item Restart $M$ independent paths at $x_0^\varepsilon$, sharing $K_{\varepsilon,k}$ and
$\theta_{\varepsilon,k}$. Hold these parameters fixed while sampling $J$ one-step transitions per
path under \eqref{eq:rl_infinite_gaussian_policy} and observing their running costs.
\item {
Use the sampled increments in \eqref{eq:rl_infinite_critic_update} to update
$\theta_{\varepsilon,k+1}$.}
\item {
On the same paths, evaluate \eqref{eq:rl_infinite_actor_increment} with the updated critic
and update $K_{\varepsilon,k+1}$ by \eqref{eq:rl_infinite_actor_update}.}
\end{enumerate}

\medskip
\textbf{Output.}
Return {the infinite-horizon policy} \eqref{eq:rl_infinite_gaussian_policy} and its nonnegative conditional
mean \eqref{eq:rl_infinite_stationary_mean}.
\end{algorithmblock}
}

{The length $J\Delta t$ is an observation window, not a control horizon; its endpoint
contributes no terminal return or continuation value. Restarting from nonzero states supplies data
before the second moment decays under a stable policy, as shown in
\eqref{eq:infinite_second_moment_solution}. A nonzero path stays in its initial region, so the two
actors require separate positive- and negative-state runs. At an iterate satisfying
\eqref{eq:rl_infinite_fixed_policy_stability}, the critic has the total-cost interpretation in
\eqref{eq:rl_infinite_fixed_policy_value}. }

\section{Numerical experiments}
\label{sec:numerical_experiments}

{We report separate finite- and infinite-horizon experiments. The finite-horizon study
compares learned coefficients with analytical targets and examines running costs and mean-policy
compatibility. The infinite-horizon study evaluates the {batch
Algorithm~\ref{alg:rl_infinite_actor_critic} against the regional algebraic targets} and checks the
second-moment rates $N_\pm^\Psi(\widehat K_\pm)$ from \eqref{eq:infinite_rates}. Its finite data
windows contribute only one-step increments, with no assigned terminal value. Gaussian sampling
implements the prescribed moments; the control problems retain the general moment-constrained policy
class.}

\subsection{Finite-horizon design}

{For the finite-horizon study, we use time-dependent dynamics based on
\citet[Section~4]{WuEtAl2020}. The horizon is $T=0.1$, with one state, three controls, and three
Brownian components. For $t\in[0,T]$,}
\begin{align}
A(t)&=0.2+t,\nonumber\\
B(t)&=(-5+t,-10+t,20+t)^\top,\nonumber\\
C(t)&=(-0.84+t,-3.78+t,0.849+t)^\top,\nonumber\\
D(t)&=
\begin{pmatrix}
6.85+t&11.22+t&-1.98+t\\
-8.78+t&13.24+t&-5.44+t\\
0.68+t&14.53+t&-2.32+t
\end{pmatrix}.
\label{eq:numerical_dynamics_direct}
\end{align}
{We examine four cost specifications. The \emph{terminal} benchmark sets}
\[
R=S=q=0,\qquad q_T=1.
\]
The \emph{running} case sets
\begin{equation}
\begin{aligned}
R(t)&=\operatorname{diag}(3+t,5+t,4+t),\\
S(t)&=(0.1+t,0.4+t,0.5+t)^\top,\\
q(t)&=10+t,\qquad q_T=0.
\end{aligned}
\label{eq:numerical_running_direct}
\end{equation}
For the \emph{mixed compatible} case, define $h(t):=B(t)+D(t)C(t)$ and take
\begin{equation}
\begin{aligned}
R(t)&=0.02D(t)D(t)^\top,\\
S(t)&=0.02h(t),\\
q(t)&=0.02h(t)^\top[D(t)D(t)^\top]^{-1}h(t)+1,
\qquad q_T=1.
\end{aligned}
\label{eq:numerical_compatible_direct}
\end{equation}
Thus Assumption~\ref{ass:alignment} holds with $\lambda=0.02$. The \emph{mixed incompatible} case uses
\begin{equation}
R=\operatorname{diag}(20,2,20),\qquad
S=(3,-1,0)^\top,\qquad q=1.95,\qquad q_T=1,
\label{eq:numerical_incompatible_direct}
\end{equation}
which does not satisfy Assumption~\ref{ass:alignment}.

{For each specification, $K_+^0$ denotes the model-based gain without exploration and
$K_+^{\Psi_0}$ the model-based optimal mean gain at the training covariance. Neither oracle is supplied
to the learner. We use $N=100$ time cells and $x_0=1$, so only the positive region is active. The
covariance is $\Psi_0=\tau I_3$, with $\tau=0.005$ for the terminal case, $0.01$ for the running and
compatible cases, and $0.02$ for the incompatible case.}

{Each experiment uses seeds $0,1,2$, batches of $512$ trajectories, $300$ critic warm-up
steps, and $4000$ actor iterations with one critic update per actor update. We set
$\alpha_\phi=0.01$, $\alpha_\theta=1$, the initial actor gain to $0.02$, and the initial critic
coefficient to $1$. One actor gain is shared within each of $20$ blocks of time cells. We average the
cell gradients within a block, clip the block gradient at $0.1$, and then project it onto the
nonnegative cone. Each critic-gradient coordinate is clipped to $[-5,5]$ before projection.\footnote{Block sharing reduces the number of actor parameters;
averaging, clipping, and projection control individual updates. These computational choices leave the
control problem's policy class unchanged.} The critic directly learns the coefficient in
\eqref{eq:rl_param_value}:}
\begin{equation}
V^\theta(t_i,x)=\theta_{+,i}x^2,\qquad x>0,
\label{eq:numerical_direct_critic}
\end{equation}

For a learned gain $\widehat K$ and a reference gain $\overline K$, define
\begin{equation}
\mathcal E(\widehat K;\overline K)
:=\frac{\left(\sum_{i=0}^{N-1}\|\widehat K(t_i)-\overline K(t_i)\|^2\Delta t\right)^{1/2}}
{\max\left\{\left(\sum_{i=0}^{N-1}\|\overline K(t_i)\|^2\Delta t\right)^{1/2},1\right\}}.
\label{eq:numerical_metrics_direct}
\end{equation}
{The errors $\mathcal E(\widehat K;K_+^{\Psi_0})$ and $\mathcal E(\widehat K;K_+^0)$ compare the actor in
\eqref{eq:rl_gaussian_policy} with the exploratory and original targets, respectively. The critic error
is the relative root-mean-square error of $\theta_{+,i}$ against the analytical coefficient defined by
\eqref{eq:rl_positive_theta_equation} for the current policy.}

The \emph{terminal-only ablation} deletes the running-cost contribution from both the critic return in
\eqref{eq:rl_realized_return} and the actor increment in \eqref{eq:rl_martingale_increment}, while
keeping all other settings unchanged. It is reported only in the running case, where $q_T=0$
isolates the omitted running contribution.

\subsection{Finite-horizon coefficient and policy comparison}

Table~\ref{tab:numerical_direct_results} reports the final errors as the mean $\pm$ sample standard
deviation over three seeds, with $\mathcal E_{\Psi_0}:=\mathcal E(\widehat K;K_+^{\Psi_0})$ and
$\mathcal E_0:=\mathcal E(\widehat K;K_+^0)$.

\begin{table}[H]
\centering
\scriptsize
\setlength{\tabcolsep}{2.3pt}
{
\begin{tabular}{l l c c c}
\hline
Case & Method & Critic error & $\mathcal E_{\Psi_0}$ & $\mathcal E_0$\\
\hline
Terminal & Full & $0.0138\pm0.0039$ & $0.0035\pm0.0004$ & $0.0035\pm0.0004$\\
Running & Full & $0.0263\pm0.0136$ & $0.0071\pm0.0006$ & $0.0070\pm0.0006$\\
Running & Ablation & $1.0000\pm6.9\times10^{-8}$ & $0.0683\pm0.0001$ & $0.0678\pm0.0001$\\
Compatible & Full & $0.1004\pm0.0626$ & $0.0049\pm0.0002$ & $0.0049\pm0.0002$\\
Incompatible & Full & $0.2855\pm0.1225$ & $0.0098\pm0.0007$ & $0.0115\pm0.0003$\\
\hline
\end{tabular}
}
\caption{{Finite-horizon errors after $4000$ actor iterations using the unnormalized gradients in \eqref{eq:rl_critic_sample_gradient} and \eqref{eq:rl_actor_gradient}. Entries are means $\pm$ sample standard deviations over three independent seeds.}}
\label{tab:numerical_direct_results}
\end{table}

{In the terminal benchmark, the two oracle gains agree to numerical precision, and the
learned gain has the same error, $0.0035$, relative to both. This benchmark has $R=0$ and therefore
does not satisfy Assumption~\ref{ass:control_weight_positive_definite}. In the running-cost case, the
full critic error is $0.0263$, compared with $1.0000$ under terminal-only ablation; the actor error to
$K_+^{\Psi_0}$ rises from $0.0071$ to $0.0683$. Because $q_T=0$, deleting the running cost makes the
ablated return identically zero. This comparison isolates the contribution of running costs in this
specification.}

{The mixed-cost cases examine the mean-policy identity in
Theorem~\ref{thm:alignment_consistency}. In the compatible case,}
\[
\|K_+^{\Psi_0}-K_+^0\|_{L^2}=1.95\times10^{-17},
\]
{and the full method has the same error, $0.0049$, relative to both oracles. The corresponding critic
error is $0.1004\pm0.0626$.} In the incompatible case,
\[
{\|K_+^{\Psi_0}-K_+^0\|_{L^2}=0.0078.}
\]
{The final gain is closer to the exploratory oracle than to
the original oracle, with errors $0.0098$ and $0.0115$; this ordering holds for all three seeds.
Together, the two cases illustrate exact oracle agreement under compatibility and a distinct
exploratory target in this incompatible specification.}

\begin{figure}[t]
\centering
\includegraphics[width=\linewidth]{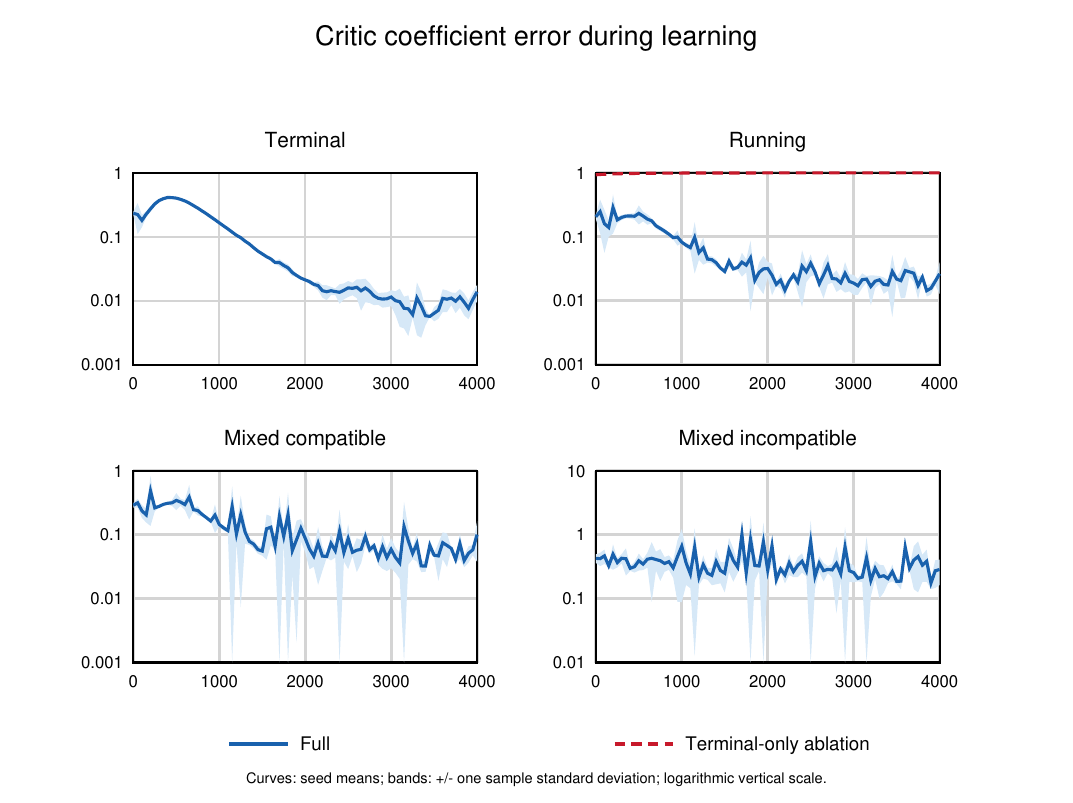}
\caption{Critic coefficient error during learning. The full method is shown by solid lines and the terminal-only ablation by dashed lines. Curves are three-seed means, shaded bands are one sample standard deviation, and the vertical scale is logarithmic.}
\label{fig:numerical_direct_critic_errors}
\end{figure}

\begin{figure}[t]
\centering
\includegraphics[width=\linewidth]{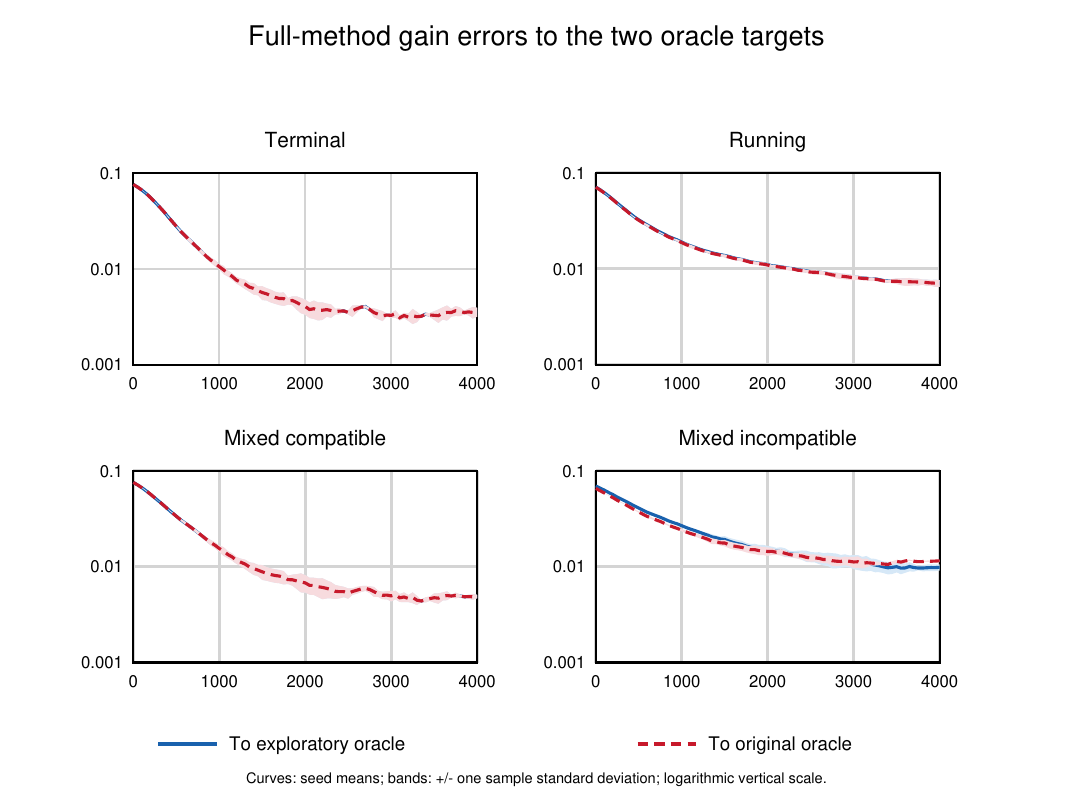}
\caption{Full-method gain errors relative to the exploratory and original oracles. The error to $K_+^{\Psi_0}$ is shown by solid lines and the error to $K_+^0$ by dashed lines. Curves are three-seed means, shaded bands are one sample standard deviation, and the vertical scale is logarithmic.}
\label{fig:numerical_direct_gain_errors}
\end{figure}

{Figure~\ref{fig:numerical_direct_critic_errors} shows that the terminal and running full-method critic
errors end at the order of $10^{-2}$, whereas the running ablation remains near $10^0$; the compatible
and incompatible curves have wider three-seed bands. In Figure~\ref{fig:numerical_direct_gain_errors},
the two oracle-error curves overlap in the terminal and compatible cases. In the incompatible case, the
solid curve crosses below the dashed curve late in training and remains lower at iteration $4000$.}

\subsection{Finite-horizon sensitivity}

{We vary one critic parameter at a time in the incompatible case. Starting from
$(\theta_0,\alpha_\theta)=(1,1)$, we try $\theta_0=0.5,2$ and $\alpha_\theta=0.5,2$ while holding
other settings fixed. Across the five settings, the critic error ranges from $0.2363$ to $0.3736$,
$\mathcal E(\widehat K;K_+^{\Psi_0})$ from $0.0093$ to $0.0118$, and
$\mathcal E(\widehat K;K_+^0)$ from $0.0112$ to $0.0126$. The actor remains closer to the exploratory
oracle in every setting, while the critic error changes appreciably with $\alpha_\theta$. These
three-seed observations concern oracle ordering over the reported parameter range.}

\subsection{{Infinite-horizon batch actor--critic experiment}}

{The infinite-horizon study uses a separate constant-coefficient scalar model with one
control and one Brownian component:}
\begin{equation}
\begin{aligned}
A&=-0.5,\quad B=-1,\quad C=0,\quad D=1,\\
R&=1,\quad S=-1,\quad q=2,\quad \Psi=\psi=0.1.
\end{aligned}
\label{eq:numerical_infinite_model}
\end{equation}
The block cost matrix is
\[
Q=\begin{pmatrix}R&S\\S&q\end{pmatrix}
=\begin{pmatrix}1&-1\\-1&2\end{pmatrix}\succ0.
\]
Moreover, the time-independent compatibility condition \eqref{eq:infinite_alignment} holds with
$\lambda=1$, because $R=DD^\top$ and $S=B+DC=-1$. Hence Proposition~\ref{prop:infinite_alignment} applies
whenever the hypotheses of Theorem~\ref{thm:infinite_horizon} hold.

{We train Algorithm~\ref{alg:rl_infinite_actor_critic} directly on this model. The
positive and negative branches start separately from $x_0^+=1$ and $x_0^-=-1$; within each region,
{all paths in the same training round share a common gain and critic coefficient}. The sampling policy is}
\[
U\mid X=x\sim
\begin{cases}
\mathcal N(xK_+,x^2\psi),&x>0,\\
\mathcal N((-x)K_-,x^2\psi),&x<0.
\end{cases}
\]
As in \eqref{eq:rl_gaussian_policy}, this distribution is used to generate samples with the required
conditional moments; the theoretical policy class in Theorem~\ref{thm:infinite_horizon}
remains distribution free beyond those moments.

{For each of $20$ independent seeds, training runs for $400$ rounds. Every round restarts
$M=512$ paths per region and observes $J=50$ steps of length $\Delta t=0.01$; the resulting window
$J\Delta t=0.5$ is not a control horizon. We set $K_{\pm,0}=0.02$, $\theta_{\pm,0}=1$,
$\alpha_K=0.01$, and $\alpha_\theta=0.2$. The critic updates first, followed by the actor on the
same batch. We use no critic warm-up, state normalization, gradient clipping, or model-based stability
gate. Within each numerical cell, we hold $U_i/|X_i|$ fixed, propagate the resulting geometric state
equation, and approximate running cost at the left endpoint. Neither update uses a terminal value or
finite-window return.}

For \eqref{eq:numerical_infinite_model}, the quantities in \eqref{eq:infinite_rates} reduce to
\[
\begin{aligned}
N_\pm^\Psi(K)
&=2A+C^2+D^2\psi+D^2K^2\pm2(B+DC)K,\\
L_\pm^\Psi(K)
&=q+R\psi+RK^2\pm2SK.
\end{aligned}
\]
{{The positive-region infinite-horizon oracle is}
$(K_+^{\Psi,*},P_+^{\Psi,*})=(1.0000,0.5789)$. {The negative-region oracle is}
$(K_-^{\Psi,*},P_-^{\Psi,*})=(0.0000,2.3333)$. Their second-moment rates are $-1.9000$ and $-0.9000$.
These model-based quantities are used only to assess the learned policy after training.
Table~\ref{tab:numerical_infinite_results} reports the outcomes.}

\begin{table}[H]

\centering
\scriptsize
\setlength{\tabcolsep}{2.7pt}
\begin{tabular}{c c c c c c c}
\hline
Region & Learned $K_\pm$ & $K_\pm^{\Psi,*}$ & Learned $\theta_\pm$ & $P_\pm^{\Psi,*}$
& \shortstack{Current-policy\\critic error} & \shortstack{Stable\\runs}\\
\hline
$+$ & $0.9837\pm0.0097$ & $1.0000$ & $0.5529\pm0.0589$ & $0.5789$ & $0.0517$ & $20/20$\\
$-$ & $0.0000\pm0.0000$ & $0.0000$ & $2.3526\pm0.0282$ & $2.3333$ & $0.0291$ & $20/20$\\
\hline
\end{tabular}
\caption{Batch infinite-horizon actor--critic results after $400$ rounds. Learned quantities are means
$\pm$ sample standard deviations over $20$ independent seeds. The current-policy critic error is the
mean of $|\widehat\theta_\pm+L_\pm^\Psi(\widehat K_\pm)/N_\pm^\Psi(\widehat K_\pm)|$; ``stable runs''
counts runs with $N_\pm^\Psi(K_{\pm,k})<0$ at every recorded iterate. Rates and oracle values are
retrospective model-based diagnostics, not training gates.}
\label{tab:numerical_infinite_results}
\end{table}

{The learned positive-region gain approaches its oracle value $1.0000$, while the
negative-region gain reaches the nonnegative boundary $0.0000$. Mean absolute critic error relative
to the current-policy value is $0.0517$ and $0.0291$ in the two regions. All $20$ runs completed
$400$ rounds, and every recorded gain had a negative diagnostic rate. This observation does not imply
that the ungated update preserves stability in general.}

\begin{figure}[H]

\centering
\includegraphics[width=0.92\linewidth]{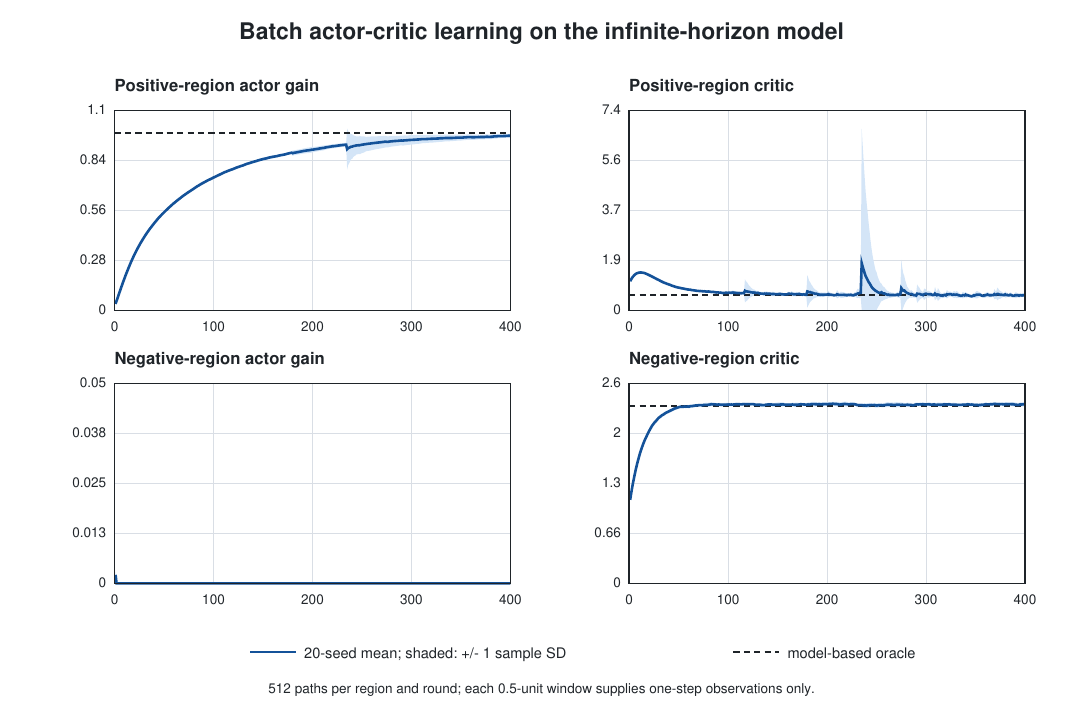}
\caption{Batch infinite-horizon learning across $20$ independent seeds. Solid curves show means, shaded
bands one sample standard deviation, and dashed lines the {model-based algebraic targets}. The positive
critic has intermittent excursions even though the final learned gains are near their targets; the
curves do not establish monotone or general convergence.}
\label{fig:numerical_infinite_batch_learning}
\end{figure}

{
To compare the learning variability and recorded policy stability of single-path and
multiple-path sampling, we repeat the experiment with one newly initialized path per region and
round. The dynamics, initial parameters, $400$ rounds, observation window, and learning rates are
unchanged. All $20$ runs complete training. The final gains are $0.8114\pm0.2761$ and
$0.0287\pm0.0620$; the numbers of runs with all recorded second-moment rates negative are $19/20$
and $18/20$, respectively. Compared with the $512$-path results in
Table~\ref{tab:numerical_infinite_results}, the final gains have larger variability across seeds,
and fewer runs maintain negative rates at every recorded iterate.}

{The studies use different dynamics and assess different claims. The finite-horizon
study examines learning and the exploratory-versus-original mean comparison; the infinite-horizon
study assesses the {batch infinite-horizon algorithm} on a separate scalar model.}

{
\section{Conclusion}
\label{sec:conclusion}

This paper studied continuous-time stochastic LQ control with a scalar state and nonnegative
control constraints. We introduced an auxiliary problem without entropy, using randomized
policies with nonnegative conditional means and prescribed state-dependent covariance. Moment
reduction and state separation yielded two regional Riccati equations and characterized the
finite-horizon value function and optimal conditional mean (Theorem~\ref{thm:main}). Quadratic
running costs can make exploration change the optimal mean. We gave a sufficient control--cost
compatibility condition under which that mean coincides with the original optimal constrained
feedback (Theorem~\ref{thm:alignment_consistency}). For constant coefficients, we characterized
the undiscounted infinite-horizon optimum and established mean-square asymptotic stability under
the stated assumptions (Theorem~\ref{thm:infinite_horizon}).

{
For unknown dynamics, we developed finite- and infinite-horizon actor--critic schemes using the
martingale methods of \citet{JiaZhou2022PE,JiaZhou2022PG}. In the reported running-cost comparison,
including the running cost reduced both critic and policy errors. The compatible and incompatible
cases illustrated, respectively, agreement of the two optimal mean feedbacks and a distinct
exploratory target. The learned gains were close to their analytical targets in the reported finite-
and infinite-horizon experiments (Tables~\ref{tab:numerical_direct_results}
and~\ref{tab:numerical_infinite_results}).}
}

\clearpage

\appendix
\section{Proofs for the exploratory-control results}
\label{app:main_proofs}

{This appendix gives the moment-reduction, state-separation, regional-equation, and
viscosity-verification calculations used in Section~\ref{sec:main_results}.}

\subsection{{Characteristics under $\pi$}}
\label{app:exploratory_characteristics}

{{We derive the coefficients of the exploratory dynamics} in
\eqref{eq:exploratory_drift}--\eqref{eq:exploratory_state} from the prescribed covariance and the
point-action dynamics. The calculation averages the drift and quadratic variation under a
policy.}

Fix $(s,x)$ and write $a=a^\pi(s,x)$ and
$\Sigma=\Sigma^\pi(s,x)=x^2\Psi(s)$. By linearity of integration,
\begin{align}
\widetilde b^\pi(s,x)
&=\int_{\mathbb R^n}[A(s)x+B(s)^\top u] \,\pi_s(du\mid x)\nonumber\\
&=A(s)x+B(s)^\top a.
\label{eq:app_exploratory_drift_calc}
\end{align}
For the $j$th diffusion channel, write $d_j=D_{\cdot j}(s)\in\mathbb R^n$ and $c_j=C_j(s)$. Starting from
the {definition of the quadratic variation under $\pi$},
\begin{align*}
\bigl(\widetilde\sigma_j^\pi(s,x)\bigr)^2
&=\int_{\mathbb R^n}(xc_j+d_j^\top u)^2\,\pi_s(du\mid x)\\
&=x^2c_j^2
  +2xc_j d_j^\top\int_{\mathbb R^n}u\,\pi_s(du\mid x)
  +\int_{\mathbb R^n}u^\top d_jd_j^\top u\,\pi_s(du\mid x)\\
&=x^2c_j^2+2xc_jd_j^\top a
  +\int_{\mathbb R^n}u^\top d_jd_j^\top u\,\pi_s(du\mid x).
\end{align*}
The remaining quadratic term is evaluated from the covariance identity
\[
\int_{\mathbb R^n}uu^\top\,\pi_s(du\mid x)=aa^\top+\Sigma,
\qquad \Sigma=x^2\Psi(s).
\]
Indeed, using $u^\top d_jd_j^\top u=\operatorname{tr}(d_jd_j^\top uu^\top)$,
\begin{align*}
\int_{\mathbb R^n}u^\top d_jd_j^\top u\,\pi_s(du\mid x)
&=\operatorname{tr}\!\left[d_jd_j^\top
  \int_{\mathbb R^n}uu^\top\,\pi_s(du\mid x)\right]\\
&=\operatorname{tr}[d_jd_j^\top(aa^\top+\Sigma)]\\
&=(d_j^\top a)^2+d_j^\top\Sigma d_j.
\end{align*}
Substituting this identity and then using $\Sigma=x^2\Psi(s)$ gives
\begin{align}
\bigl(\widetilde\sigma_j^\pi(s,x)\bigr)^2
&=x^2c_j^2+2xc_jd_j^\top a+(d_j^\top a)^2+d_j^\top\Sigma d_j\nonumber\\
&=(xc_j+d_j^\top a)^2+d_j^\top\Sigma d_j\nonumber\\
&=[xC_j(s)+D_{\cdot j}(s)^\top a]^2
  +x^2D_{\cdot j}(s)^\top\Psi(s)D_{\cdot j}(s).
\label{eq:app_exploratory_diffusion_calc}
\end{align}
{The first term evaluates diffusion at the policy mean; the second is the variance from
action randomization. Summing over Brownian components gives}
\begin{align}
\sum_{j=1}^m\bigl(\widetilde\sigma_j^\pi(s,x)\bigr)^2
={}&\|xC(s)+D(s)^\top a\|^2
+x^2\operatorname{tr}[D(s)D(s)^\top\Psi(s)].
\label{eq:app_exploratory_total_variance_calc}
\end{align}
{Each $\widetilde\sigma_j^\pi$ in \eqref{eq:exploratory_state} can therefore be taken as
the nonnegative square root in \eqref{eq:app_exploratory_diffusion_calc}. The construction averages
the quadratic variation $\sigma_j^2$, rather than the signed diffusion coefficient.}

\subsection{Moment reduction and state separation}

We give the details of Proposition~\ref{prop:distribution_to_mean}. Fix $(t,x)$ and an admissible law
$\pi\in\mathcal P_{\Psi,+}(t,x)$, and write
\[
a:=\int_{\mathbb R^n}u\,\pi(du),
\qquad
\Sigma:=\operatorname{Cov}_\pi(u)=x^2\Psi(t).
\]
By the definition of covariance,
\begin{align*}
\Sigma
&=\int_{\mathbb R^n}(u-a)(u-a)^\top\pi(du)\\
&=\int_{\mathbb R^n}uu^\top\pi(du)
-a\int_{\mathbb R^n}u^\top\pi(du)
-\left(\int_{\mathbb R^n}u\,\pi(du)\right)a^\top
+aa^\top\\
&=\int_{\mathbb R^n}uu^\top\pi(du)-aa^\top.
\end{align*}
Hence
\begin{equation}
\int_{\mathbb R^n}uu^\top\pi(du)=aa^\top+\Sigma=aa^\top+x^2\Psi(t),
\label{eq:app_second_moment_detailed}
\end{equation}
which proves \eqref{eq:second_moment_reduction}.

The drift reduction is immediate from the affine dependence on the action:
\begin{align*}
\int_{\mathbb R^n}b(t,x,u)\pi(du)
&=\int_{\mathbb R^n}[A(t)x+B(t)^\top u]\pi(du)\\
&=A(t)x+B(t)^\top a,
\end{align*}
which is \eqref{eq:drift_reduction}.

For the diffusion term, first note that
\begin{align*}
\sum_{j=1}^m\sigma_j^2(t,x,u)
&=\|xC(t)+D(t)^\top u\|^2\\
&=x^2\|C(t)\|^2+2xC(t)^\top D(t)^\top u
+u^\top D(t)D(t)^\top u.
\end{align*}
Taking expectation under $\pi$ and using \eqref{eq:app_second_moment_detailed},
\begin{align*}
\int u^\top D D^\top u\,\pi(du)
&=\operatorname{tr}\!\left(DD^\top\int uu^\top\pi(du)\right)\\
&=\operatorname{tr}[DD^\top(aa^\top+\Sigma)]\\
&=a^\top DD^\top a+\operatorname{tr}(DD^\top\Sigma).
\end{align*}
Therefore,
\begin{align}
\sum_{j=1}^m\int\sigma_j^2(t,x,u)\pi(du)
&=x^2\|C\|^2+2xC^\top D^\top a+a^\top DD^\top a
+\operatorname{tr}(DD^\top\Sigma)\nonumber\\
&=\|xC+D^\top a\|^2
+x^2\operatorname{tr}[DD^\top\Psi(t)],
\label{eq:variance_reduction_calculation}
\end{align}
where the time argument $t$ has been suppressed in the second line only for readability.
This proves \eqref{eq:variance_reduction}.

For the running cost, the cross and state terms depend only on the first moment, while the control
term is handled by the same trace identity:
\begin{align*}
\int u^\top R(t)u\,\pi(du)
&=\operatorname{tr}\!\left(R(t)\int uu^\top\pi(du)\right)\\
&=a^\top R(t)a+\operatorname{tr}[R(t)\Sigma]\\
&=a^\top R(t)a+x^2\operatorname{tr}[R(t)\Psi(t)].
\end{align*}
Consequently,
\begin{align}
\int\ell(t,x,u)\pi(du)
={}&a^\top R(t)a+2xS(t)^\top a+q(t)x^2
+x^2\operatorname{tr}[R(t)\Psi(t)],
\label{eq:cost_reduction_calculation}
\end{align}
which is \eqref{eq:cost_reduction}.

{Under the prescribed covariance $x^2\Psi(t)$, the HJB integrand depends on a law only
through its mean $a$. Every $a\in\mathbb R_+^n$ is realized by, for example,
$\mathcal N(a,x^2\Psi(t))$. This Gaussian law establishes attainability without restricting the
admissible policy class. Hence the infimum over laws in \eqref{eq:distributional_hjb} equals the
infimum over means in \eqref{eq:mean_hjb}, proving Proposition~\ref{prop:distribution_to_mean}.}

{For state separation, fix $p\ge0$.} If $x>0$, set $a=xK$. Since $\mathbb R_+^n$ is a cone, $a\in\mathbb
R_+^n$ if and only if $K\in\mathbb R_+^n$, and
\[
a^\top\widehat R(t,p)a+2x\widehat S(t,p)^\top a
=x^2\{K^\top\widehat R(t,p)K+2\widehat S(t,p)^\top K\}.
\]
If $x<0$, set $a=(-x)K$ with $-x>0$. Then
\[
a^\top\widehat R(t,p)a+2x\widehat S(t,p)^\top a
=x^2\{K^\top\widehat R(t,p)K-2\widehat S(t,p)^\top K\}.
\]
Finally, if $x=0$, the objective reduces to $a^\top\widehat R(t,p)a$.
{Assumption~\ref{ass:control_weight_positive_definite} gives $\widehat R(t,p)=R(t)+pD(t)D(t)^\top\succ0$
for every $p\ge0$. Hence its unique minimum over $\mathbb R_+^n$ is attained at $a=0$, and the minimizer
in each of the three state regions is unique.} These three cases prove
\eqref{eq:state_separation_value}--\eqref{eq:state_separation_minimizer}.

\subsection{Derivation of the regional Riccati equations}

For $x>0$, take the regional candidate $v(t,x)=P_+^\Psi(t)x^2$. Then
\[
v_t=\dot P_+^\Psi x^2,\qquad v_x=2P_+^\Psi x,\qquad v_{xx}=2P_+^\Psi.
\]
Substituting these derivatives into \eqref{eq:mean_hjb} and applying \eqref{eq:state_separation_value}
gives the $+$ equation in \eqref{eq:riccati_pm}. For $x<0$, the same calculation with
$v(t,x)=P_-^\Psi(t)x^2$ and the negative branch of \eqref{eq:state_separation_value} gives the $-$
equation. The terminal condition follows from $v(T,x)=q_Tx^2$.

\subsection{Uniqueness of the regional minimizers}

{For $p\ge0$, Assumption~\ref{ass:control_weight_positive_definite} gives
\[
\widehat R(t,p)=R(t)+pD(t)D(t)^\top\succ0.
\]
Each regional objective in \eqref{eq:regional_gain} is strictly convex and coercive on the closed
convex set $\mathbb R_+^n$ and thus has a unique minimizer $K_\pm(t,p)$.} {Uniqueness of
the HJB solution is separately covered by the comparison conditions in Theorem~\ref{thm:main}.}

\subsection{Verification of Theorem~\ref{thm:main}}

{Let $\widetilde V^\Psi$ be the piecewise-quadratic candidate in
\eqref{eq:main_value}. It satisfies the terminal condition. We verify the HJB equation on each open
region and then at their common boundary.}

\textit{1. Interior regions.} For $x>0$,
\begin{equation}
(\widetilde V^\Psi)_t=\dot P_+^\Psi x^2,
\qquad
(\widetilde V^\Psi)_x=2P_+^\Psi x,
\qquad
(\widetilde V^\Psi)_{xx}=2P_+^\Psi.
\label{eq:plus_derivatives}
\end{equation}
Denote the left-hand side of \eqref{eq:mean_hjb} by $\mathrm{LHS}$. Substitution yields
\begin{align}
\mathrm{LHS}={}&x^2\big\{\dot P_+^\Psi+[2A+\|C\|^2]P_+^\Psi+q+\operatorname{tr}(R\Psi)
+P_+^\Psi\operatorname{tr}(DD^\top\Psi)\big\}
\nonumber\\
&+\min_{a\in\mathbb R_+^n}
\big\{a^\top\widehat R(t,P_+^\Psi)a
+2x\widehat S(t,P_+^\Psi)^\top a\big\}.
\label{eq:plus_hjb_substitution}
\end{align}
By \eqref{eq:state_separation_value}, the minimum equals $x^2F_+(t,P_+^\Psi)$; hence
\eqref{eq:riccati_pm} gives $\mathrm{LHS}=0$. The minimizing mean is $xK_+^\Psi(t)$. For $x<0$, the same
calculation with the negative regional equation gives $\mathrm{LHS}=0$ and the minimizing mean
$-xK_-^\Psi(t)$; the details are identical and are omitted.

\textit{2. Separating boundary.} {At $x=0$, the candidate is continuous because}
\[
\widetilde V^\Psi(t,0)=P_+^\Psi(t)0^2=P_-^\Psi(t)0^2=0.
\]
Its first-order derivatives also exist and satisfy
\[
(\widetilde V^\Psi)_t(t,0)=0,
\qquad
(\widetilde V^\Psi)_x(t,0)=0.
\]
{When $P_+^\Psi(t)\neq P_-^\Psi(t)$, the second derivative may jump at the boundary.
We use the viscosity verification of \citet{WuEtAl2020}. Let $D_{t,x}^{1,2,+}$ and
$D_{t,x}^{1,2,-}$ denote the parabolic second-order superdifferential and subdifferential;
$(\eta,p,P)$ records their time, first-state, and second-state components. At $x=0$, we obtain}
\begin{equation}
\begin{aligned}
D_{t,x}^{1,2,+}\widetilde V^\Psi(t,0)
&=\{0\}\times\{0\}\times[P_+^\Psi(t),\infty),\\
D_{t,x}^{1,2,-}\widetilde V^\Psi(t,0)
&=\{0\}\times\{0\}\times(-\infty,P_-^\Psi(t)].
\end{aligned}
\label{eq:boundary_semijets}
\end{equation}

Define the HJB term
\begin{align}
W^\Psi(t,x,a,p,P):={}&\frac12P\|xC(t)+D(t)^\top a\|^2
\nonumber\\
&+\frac12Px^2\operatorname{tr}[D(t)D(t)^\top\Psi(t)]
\nonumber\\
&+p[A(t)x+B(t)^\top a]
+a^\top R(t)a+2xS(t)^\top a
\nonumber\\
&+q(t)x^2
+x^2\operatorname{tr}[R(t)\Psi(t)].
\label{eq:hjb_operator}
\end{align}
At $x=0$, all terms containing $x$ vanish. Hence, for any
$(\eta,p,P)\in D_{t,x}^{1,2,+}\widetilde V^\Psi(t,0)$,
\begin{align}
\eta+\min_{a\in\mathbb R_+^n}W^\Psi(t,0,a,p,P)
&=\min_{a\in\mathbb R_+^n}
\frac12a^\top\bigl[P D(t)D(t)^\top+2R(t)\bigr]a
\nonumber\\
&\ge
\min_{a\in\mathbb R_+^n}
\frac12a^\top\bigl[P_+^\Psi(t)D(t)D(t)^\top+2R(t)\bigr]a
=0,
\label{eq:boundary_subsolution}
\end{align}
where the last equality follows by taking $a=0$ and from the positive semidefiniteness of the quadratic
form. Thus $\widetilde V^\Psi$ is a viscosity subsolution on the separating boundary.

On the other hand, for any $(\eta,p,P)\in D_{t,x}^{1,2,-}\widetilde V^\Psi(t,0)$,
\begin{align}
\eta+\min_{a\in\mathbb R_+^n}W^\Psi(t,0,a,p,P)
&=\min_{a\in\mathbb R_+^n}
\frac12a^\top\bigl[P D(t)D(t)^\top+2R(t)\bigr]a
\nonumber\\
&\le
\min_{a\in\mathbb R_+^n}
\frac12a^\top\bigl[P_-^\Psi(t)D(t)D(t)^\top+2R(t)\bigr]a
=0.
\label{eq:boundary_supersolution}
\end{align}
Hence $\widetilde V^\Psi$ is also a viscosity supersolution on the separating boundary. Together with the
verification in the two interior regions, this shows that $\widetilde V^\Psi$ is a viscosity solution of
\eqref{eq:mean_hjb} with the terminal condition $V(T,x)=q_Tx^2$.

Finally, for $(t,0)$, take
\[
(\eta^*(t,0),p^*(t,0),P^*(t,0))=(0,0,P_+^\Psi(t))
\in D_{t,x}^{1,2,+}\widetilde V^\Psi(t,0)
\]
and $a^{\Psi,*}(t,0)=0$. Then
\[
\eta^*(t,0)+W^\Psi\bigl(t,0,a^{\Psi,*}(t,0),p^*(t,0),P^*(t,0)\bigr)=0.
\]
The regional minimizing means together with $a^{\Psi,*}(t,0)=0$ yield \eqref{eq:optimal_mean}. Under the
comparison condition stated in Theorem~\ref{thm:main}, the viscosity solution is
unique, and hence $\widetilde V^\Psi=V^\Psi$.

\section{Derivations for the unknown-dynamics learning section}
\label{app:rl_derivations}

{This appendix derives the fixed-policy regional equations used in
Section~\ref{sec:unknown_dynamics}.}

\subsection{Fixed-policy regional equations}

Consider first $x>0$, for which $a^\phi(t,x)=x\phi_+(t)$.
Proposition~\ref{prop:distribution_to_mean} gives
\begin{equation}
A(t)x+B(t)^\top a^\phi(t,x)
=x\left[A(t)+B(t)^\top\phi_+(t)\right],
\label{eq:app_rl_positive_drift}
\end{equation}
\begin{align}
&\left\|xC(t)+D(t)^\top a^\phi(t,x)\right\|^2
+x^2\operatorname{tr}\!\left[D(t)D(t)^\top\Psi_0\right]
\nonumber\\
&\qquad=
x^2\left\{
\left\|C(t)+D(t)^\top\phi_+(t)\right\|^2
+\operatorname{tr}\!\left[D(t)D(t)^\top\Psi_0\right]
\right\},
\label{eq:app_rl_positive_variance}
\end{align}
and
\begin{align}
\overline\ell^\phi(t,x)
=x^2\bigg\{
&\phi_+(t)^\top R(t)\phi_+(t)
+2S(t)^\top\phi_+(t)+q(t)
\nonumber\\
&+\operatorname{tr}\!\left[R(t)\Psi_0\right]
\bigg\}.
\label{eq:app_rl_positive_cost}
\end{align}
Substituting $V(t,x)=\theta_+^\phi(t)x^2$ and
\eqref{eq:app_rl_positive_drift}--\eqref{eq:app_rl_positive_cost} into \eqref{eq:rl_fixed_policy_pde}
gives \eqref{eq:rl_positive_theta_equation}.  For $x<0$, the same calculation with
$a^\phi(t,x)=-x\phi_-(t)$ and $V(t,x)=\theta_-^\phi(t)x^2$ gives \eqref{eq:rl_negative_theta_equation}.
At $x=0$, the policy is $\delta_0$, the state and all costs remain zero.  This verifies
\eqref{eq:rl_candidate_fixed_policy_value}.  Using \eqref{eq:regional_matrices}, the control-dependent
terms in the two coefficient equations are
\begin{equation}
\phi_+(t)^\top\widehat R\!\left(t,\theta_+^\phi(t)\right)\phi_+(t)
+2\widehat S\!\left(t,\theta_+^\phi(t)\right)^\top\phi_+(t)
\label{eq:app_rl_plus_bridge}
\end{equation}
and
\begin{equation}
\phi_-(t)^\top\widehat R\!\left(t,\theta_-^\phi(t)\right)\phi_-(t)
-2\widehat S\!\left(t,\theta_-^\phi(t)\right)^\top\phi_-(t),
\label{eq:app_rl_minus_bridge}
\end{equation}
respectively. Minimization over $\mathbb R_+^n$ therefore recovers the regional quadratic programs
{\eqref{eq:regional_value}--\eqref{eq:regional_gain}}.

\end{document}